\documentclass{gtart}

\def\ifplaintex{\expandafter\ifx\csname documentclass\endcsname\relax}

\def\gtp{{\mathsurround=0pt\it $\cal G\mskip-2mu$eometry \&\ 
$\cal T\!\!$opology $\cal P\!$ublications}}  

\def\recd{{\small Received:\qua\receiveddate\ifx\reviseddate\relax
\else\qquad Revised:\qua\reviseddate\fi\par}} 

\def\lognumber#1{\def\thelognumber{#1}}
\def\volumenumber#1{\def\thevolumenumber{#1}}
\def\volumeyear#1{\def\thevolumeyear{#1}}
\def\papernumber#1{\def\thepapernumber{#1}}
\def\pagenumbers#1#2{\def\startpage{#1}\def\finishpage{#2}}
\def\published#1{\def\publishdate{#1}}

\def\received#1{\def\receiveddate{#1}}

\def\accepted#1{\def\accepteddate{#1}}
\def\asciititle#1{\def\theasciititle{#1}}
\def\covertitle#1{\def\thecovertitle{#1}}

\def\asciiaddress#1{\def\theasciiaddress{#1}}

\long\def\asciiabstract#1{\long\def\theasciiabstract{#1}}
\def\asciikeywords#1{\def\theasciikeywords{#1}}

\let\\\par\let\thelognumber\relax\let\thevolumenumber\relax
\let\thepapernumber\relax\let\thevolumeyear\relax\let\startpage\relax
\let\finishpage\relax\let\publishdate\relax\let\receiveddate\relax
\let\reviseddate\relax\let\accepteddate\relax\let\theasciititle\relax
\let\thecovertitle\relax\let\theasciiauthors\relax\let\theasciiaddress\relax
\let\theasciiabstract\relax\let\theasciikeywords\relax

\let\theasciiemail\relax

\ifplaintex
\font\logobig=cmssbx10 scaled 3836
\font\logomed=cmssbx10 scaled 2557
\else
\font\logobig=cmssbx10 scaled 4200
\font\logomed=cmssbx10 scaled 2800
\fi

\long\def\makeagttitle{   
\count0=\startpage
\agt\hfill      
\hbox to 45truept{\vbox to 0pt{\vglue -13truept{\logomed A\kern -.37em{\logobig 
T}\kern -.38em G}\vss}\hss}
\break
{\small Volume \thevolumenumber\ (\thevolumeyear)
\startpage--\finishpage\nl
Published: \publishdate}

\vglue .25truein

{\parskip=0pt\leftskip 0pt plus
1fil\def\\{\par\smallskip}{\Large\bf\thetitle}\par\medskip} \vglue
0.05truein

{\parskip=0pt\leftskip 0pt plus 1fil\def\\{\par}{\sc\theauthors}
\par\medskip}%
 
\vglue 0.03truein 

{\small\leftskip 25truept\rightskip 25truept{\bf Abstract}\stdspace\theabstract

{\bf AMS Classification}\stdspace\theprimaryclass
\ifx\thesecondaryclass\relax\else; \thesecondaryclass\fi\par
{\bf Keywords}\stdspace \thekeywords\par}\vglue 7truept

}   

\ifplaintex
\font\phead=cmsl9 scaled 950
\font\pnum=cmbx10 scaled 913
\font\pfoot=cmsl9 scaled 950
\headline{\vbox to 0pt{\vskip -4.5mm\line{\small\phead\ifnum
\count0=\startpage ISSN 1472-2739 (on-line) 1472-2747 (printed)
\hfill {\pnum\folio}\else\ifodd\count0\def\\{ }%
\ifx\theshorttitle\relax\thetitle\else\theshorttitle\fi\hfill{\pnum\folio}
\else\def\\{ and }{\pnum\folio}\hfill\ifx\theshortauthors\relax\theauthors
\else\theshortauthors\fi\fi\fi}\vss}}
\footline{\vbox to 0pt{\vglue 0mm\line{\small\pfoot\ifnum\count0=\startpage
\copyright\ \gtp\hfill\else
\agt, Volume \thevolumenumber\ (\thevolumeyear)\hfill\fi}\vss}}
\else
\font=cmsl9 scaled 1050
\font\lnum=cmbx10 
\font=cmsl9 scaled 1050
\makeatletter
\def\@oddhead{{\small\lhead\ifnum\count0=\startpage ISSN 1472-2739 
(on-line) 1472-2747 (printed)\hfill {\lnum\number\count0}\else\ifodd\count0
\def\\{ }\ifx\theshorttitle\relax \thetitle \else\theshorttitle\fi\hfill
{\lnum\number\count0}\else\def\\{ and }{\lnum\number\count0}
\hfill\ifx\theshortauthors\relax 
\theauthors\else\theshortauthors\fi\fi\fi}}\def\@evenhead{\@oddhead}
\def\@oddfoot{\small\lfoot\ifnum\count0=\startpage\copyright\ \gtp\hfill\else
\agt, Volume \thevolumenumber\ (\thevolumeyear)\hfill\fi}
\def\@evenfoot{\@oddfoot}
\makeatother
\fi
\let\maketitlepage\makeagttitle
\let\makeshorttitle\maketitlepage
\let\maketitle\maketitlepage

\newwrite\gtoutfile
\long\gdef\makeheadfile{  
{\def\\{, }\def\s{ }
\immediate\openout\gtoutfile head.xxx
\immediate\write\gtoutfile{To: math@arxiv.org}
\immediate\write\gtoutfile{Subject: put OR rep NNNNN:ppppp}
\immediate\write\gtoutfile{--text follows this line--}
\immediate\write\gtoutfile{Proxy-for: \ifx\theasciiauthors\relax
\theauthors\else\theasciiauthors\fi\s<\ifx\theasciiemail\relax\theemail\else\theasciiemail\fi>}
\immediate\write\gtoutfile{\noexpand\\}
\immediate\write\gtoutfile{Authors: \ifx\theasciiauthors\relax
\theauthors\else\theasciiauthors\fi}
{\def\\{ }\immediate\write\gtoutfile{Title: \ifx\theasciititle\relax
\thetitle\else\theasciititle\fi}}
\immediate\write\gtoutfile{Subj-class: GT or SG, GR etc}
\immediate\write\gtoutfile{MSC-class: \theprimaryclass\ifx\thesecondaryclass\relax\else, \thesecondaryclass\fi}
\immediate\write\gtoutfile{Journal-ref: Algebr. Geom. Topol. \thevolumenumber\s
(\thevolumeyear) \startpage-\finishpage}
\immediate\write\gtoutfile{Comments: Published by Algebraic and
Geometric Topology at}
\immediate\write\gtoutfile{\s\s\s  http://www.maths.warwick.ac.uk/agt/AGTVol\thevolumenumber/agt-\thevolumenumber-\thepapernumber.abs.html}
\immediate\write\gtoutfile{\noexpand\\}
\immediate\write\gtoutfile{}
\ifx\theasciiabstract\relax
\immediate\write\gtoutfile{\theabstract}\else
\immediate\write\gtoutfile{\theasciiabstract}\fi
\immediate\write\gtoutfile{}
\immediate\write\gtoutfile{\noexpand\\}
\immediate\write\gtoutfile{}
\immediate\closeout\gtoutfile}}  

\def\maketitlepage{\makeagttitle\makeheadfile}
\let\makeshorttitle\maketitlepage
\let\maketitle\maketitlepage

\def\ifplaintex{\expandafter\ifx\csname documentclass\endcsname\relax}

\def\gtp{{\mathsurround=0pt\it $\cal G\mskip-2mu$eometry \&\ 
$\cal T\!\!$opology $\cal P\!$ublications}}  

\def\recd{{\small Received:\qua\receiveddate\ifx\reviseddate\relax
\else\qquad Revised:\qua\reviseddate\fi\par}} 

\def\lognumber#1{\def\thelognumber{#1}}
\def\volumenumber#1{\def\thevolumenumber{#1}}
\def\volumeyear#1{\def\thevolumeyear{#1}}
\def\papernumber#1{\def\thepapernumber{#1}}
\def\pagenumbers#1#2{\def\startpage{#1}\def\finishpage{#2}}
\def\published#1{\def\publishdate{#1}}

\def\received#1{\def\receiveddate{#1}}

\def\accepted#1{\def\accepteddate{#1}}
\def\asciititle#1{\def\theasciititle{#1}}
\def\covertitle#1{\def\thecovertitle{#1}}

\def\asciiaddress#1{\def\theasciiaddress{#1}}

\long\def\asciiabstract#1{\long\def\theasciiabstract{#1}}
\def\asciikeywords#1{\def\theasciikeywords{#1}}

\let\\\par\let\thelognumber\relax\let\thevolumenumber\relax
\let\thepapernumber\relax\let\thevolumeyear\relax\let\startpage\relax
\let\finishpage\relax\let\publishdate\relax\let\receiveddate\relax
\let\reviseddate\relax\let\accepteddate\relax\let\theasciititle\relax
\let\thecovertitle\relax\let\theasciiauthors\relax\let\theasciiaddress\relax
\let\theasciiabstract\relax\let\theasciikeywords\relax

\let\theasciiemail\relax

\ifplaintex
\font\logobig=cmssbx10 scaled 3836
\font\logomed=cmssbx10 scaled 2557
\else
\font\logobig=cmssbx10 scaled 4200
\font\logomed=cmssbx10 scaled 2800
\fi

\long\def\makeagttitle{   
\count0=\startpage
\agt\hfill      
\hbox to 45truept{\vbox to 0pt{\vglue -13truept{\logomed A\kern -.37em{\logobig 
T}\kern -.38em G}\vss}\hss}
\break
{\small Volume \thevolumenumber\ (\thevolumeyear)
\startpage--\finishpage\nl
Published: \publishdate}

\vglue .25truein

{\parskip=0pt\leftskip 0pt plus
1fil\def\\{\par\smallskip}{\Large\bf\thetitle}\par\medskip} \vglue
0.05truein

{\parskip=0pt\leftskip 0pt plus 1fil\def\\{\par}{\sc\theauthors}
\par\medskip}%
 
\vglue 0.03truein 

{\small\leftskip 25truept\rightskip 25truept{\bf Abstract}\stdspace\theabstract

{\bf AMS Classification}\stdspace\theprimaryclass
\ifx\thesecondaryclass\relax\else; \thesecondaryclass\fi\par
{\bf Keywords}\stdspace \thekeywords\par}\vglue 7truept

}   

\ifplaintex
\font\phead=cmsl9 scaled 950
\font\pnum=cmbx10 scaled 913
\font\pfoot=cmsl9 scaled 950
\headline{\vbox to 0pt{\vskip -4.5mm\line{\small\phead\ifnum
\count0=\startpage ISSN 1472-2739 (on-line) 1472-2747 (printed)
\hfill {\pnum\folio}\else\ifodd\count0\def\\{ }%
\ifx\theshorttitle\relax\thetitle\else\theshorttitle\fi\hfill{\pnum\folio}
\else\def\\{ and }{\pnum\folio}\hfill\ifx\theshortauthors\relax\theauthors
\else\theshortauthors\fi\fi\fi}\vss}}
\footline{\vbox to 0pt{\vglue 0mm\line{\small\pfoot\ifnum\count0=\startpage
\copyright\ \gtp\hfill\else
\agt, Volume \thevolumenumber\ (\thevolumeyear)\hfill\fi}\vss}}
\else
\font=cmsl9 scaled 1050
\font\lnum=cmbx10 
\font=cmsl9 scaled 1050
\makeatletter
\def\@oddhead{{\small\lhead\ifnum\count0=\startpage ISSN 1472-2739 
(on-line) 1472-2747 (printed)\hfill {\lnum\number\count0}\else\ifodd\count0
\def\\{ }\ifx\theshorttitle\relax \thetitle \else\theshorttitle\fi\hfill
{\lnum\number\count0}\else\def\\{ and }{\lnum\number\count0}
\hfill\ifx\theshortauthors\relax 
\theauthors\else\theshortauthors\fi\fi\fi}}\def\@evenhead{\@oddhead}
\def\@oddfoot{\small\lfoot\ifnum\count0=\startpage\copyright\ \gtp\hfill\else
\agt, Volume \thevolumenumber\ (\thevolumeyear)\hfill\fi}
\def\@evenfoot{\@oddfoot}
\makeatother
\fi
\let\maketitlepage\makeagttitle
\let\makeshorttitle\maketitlepage
\let\maketitle\maketitlepage

\newwrite\gtoutfile
\long\gdef\makeheadfile{  
{\def\\{, }\def\s{ }
\immediate\openout\gtoutfile head.xxx
\immediate\write\gtoutfile{To: math@arxiv.org}
\immediate\write\gtoutfile{Subject: put OR rep NNNNN:ppppp}
\immediate\write\gtoutfile{--text follows this line--}
\immediate\write\gtoutfile{Proxy-for: \ifx\theasciiauthors\relax
\theauthors\else\theasciiauthors\fi\s<\ifx\theasciiemail\relax\theemail\else\theasciiemail\fi>}
\immediate\write\gtoutfile{\noexpand\\}
\immediate\write\gtoutfile{Authors: \ifx\theasciiauthors\relax
\theauthors\else\theasciiauthors\fi}
{\def\\{ }\immediate\write\gtoutfile{Title: \ifx\theasciititle\relax
\thetitle\else\theasciititle\fi}}
\immediate\write\gtoutfile{Subj-class: GT or SG, GR etc}
\immediate\write\gtoutfile{MSC-class: \theprimaryclass\ifx\thesecondaryclass\relax\else, \thesecondaryclass\fi}
\immediate\write\gtoutfile{Journal-ref: Algebr. Geom. Topol. \thevolumenumber\s
(\thevolumeyear) \startpage-\finishpage}
\immediate\write\gtoutfile{Comments: Published by Algebraic and
Geometric Topology at}
\immediate\write\gtoutfile{\s\s\s  http://www.maths.warwick.ac.uk/agt/AGTVol\thevolumenumber/agt-\thevolumenumber-\thepapernumber.abs.html}
\immediate\write\gtoutfile{\noexpand\\}
\immediate\write\gtoutfile{}
\ifx\theasciiabstract\relax
\immediate\write\gtoutfile{\theabstract}\else
\immediate\write\gtoutfile{\theasciiabstract}\fi
\immediate\write\gtoutfile{}
\immediate\write\gtoutfile{\noexpand\\}
\immediate\write\gtoutfile{}
\immediate\closeout\gtoutfile}}  

\def\maketitlepage{\makeagttitle\makeheadfile}
\let\makeshorttitle\maketitlepage
\let\maketitle\maketitlepage

\def\ifplaintex{\expandafter\ifx\csname documentclass\endcsname\relax}

\def\gtp{{\mathsurround=0pt\it $\cal G\mskip-2mu$eometry \&\ 
$\cal T\!\!$opology $\cal P\!$ublications}}  

\def\recd{{\small Received:\qua\receiveddate\ifx\reviseddate\relax
\else\qquad Revised:\qua\reviseddate\fi\par}} 

\def\lognumber#1{\def\thelognumber{#1}}
\def\volumenumber#1{\def\thevolumenumber{#1}}
\def\volumeyear#1{\def\thevolumeyear{#1}}
\def\papernumber#1{\def\thepapernumber{#1}}
\def\pagenumbers#1#2{\def\startpage{#1}\def\finishpage{#2}}
\def\published#1{\def\publishdate{#1}}

\def\received#1{\def\receiveddate{#1}}

\def\accepted#1{\def\accepteddate{#1}}
\def\asciititle#1{\def\theasciititle{#1}}
\def\covertitle#1{\def\thecovertitle{#1}}

\def\asciiaddress#1{\def\theasciiaddress{#1}}

\long\def\asciiabstract#1{\long\def\theasciiabstract{#1}}
\def\asciikeywords#1{\def\theasciikeywords{#1}}

\let\\\par\let\thelognumber\relax\let\thevolumenumber\relax
\let\thepapernumber\relax\let\thevolumeyear\relax\let\startpage\relax
\let\finishpage\relax\let\publishdate\relax\let\receiveddate\relax
\let\reviseddate\relax\let\accepteddate\relax\let\theasciititle\relax
\let\thecovertitle\relax\let\theasciiauthors\relax\let\theasciiaddress\relax
\let\theasciiabstract\relax\let\theasciikeywords\relax

\let\theasciiemail\relax

\ifplaintex
\font\logobig=cmssbx10 scaled 3836
\font\logomed=cmssbx10 scaled 2557
\else
\font\logobig=cmssbx10 scaled 4200
\font\logomed=cmssbx10 scaled 2800
\fi

\long\def\makeagttitle{   
\count0=\startpage
\agt\hfill      
\hbox to 45truept{\vbox to 0pt{\vglue -13truept{\logomed A\kern -.37em{\logobig 
T}\kern -.38em G}\vss}\hss}
\break
{\small Volume \thevolumenumber\ (\thevolumeyear)
\startpage--\finishpage\nl
Published: \publishdate}

\vglue .25truein

{\parskip=0pt\leftskip 0pt plus
1fil\def\\{\par\smallskip}{\Large\bf\thetitle}\par\medskip} \vglue
0.05truein

{\parskip=0pt\leftskip 0pt plus 1fil\def\\{\par}{\sc\theauthors}
\par\medskip}%
 
\vglue 0.03truein 

{\small\leftskip 25truept\rightskip 25truept{\bf Abstract}\stdspace\theabstract

{\bf AMS Classification}\stdspace\theprimaryclass
\ifx\thesecondaryclass\relax\else; \thesecondaryclass\fi\par
{\bf Keywords}\stdspace \thekeywords\par}\vglue 7truept

}   

\ifplaintex
\font\phead=cmsl9 scaled 950
\font\pnum=cmbx10 scaled 913
\font\pfoot=cmsl9 scaled 950
\headline{\vbox to 0pt{\vskip -4.5mm\line{\small\phead\ifnum
\count0=\startpage ISSN 1472-2739 (on-line) 1472-2747 (printed)
\hfill {\pnum\folio}\else\ifodd\count0\def\\{ }%
\ifx\theshorttitle\relax\thetitle\else\theshorttitle\fi\hfill{\pnum\folio}
\else\def\\{ and }{\pnum\folio}\hfill\ifx\theshortauthors\relax\theauthors
\else\theshortauthors\fi\fi\fi}\vss}}
\footline{\vbox to 0pt{\vglue 0mm\line{\small\pfoot\ifnum\count0=\startpage
\copyright\ \gtp\hfill\else
\agt, Volume \thevolumenumber\ (\thevolumeyear)\hfill\fi}\vss}}
\else
\font=cmsl9 scaled 1050
\font\lnum=cmbx10 
\font=cmsl9 scaled 1050
\makeatletter
\def\@oddhead{{\small\lhead\ifnum\count0=\startpage ISSN 1472-2739 
(on-line) 1472-2747 (printed)\hfill {\lnum\number\count0}\else\ifodd\count0
\def\\{ }\ifx\theshorttitle\relax \thetitle \else\theshorttitle\fi\hfill
{\lnum\number\count0}\else\def\\{ and }{\lnum\number\count0}
\hfill\ifx\theshortauthors\relax 
\theauthors\else\theshortauthors\fi\fi\fi}}\def\@evenhead{\@oddhead}
\def\@oddfoot{\small\lfoot\ifnum\count0=\startpage\copyright\ \gtp\hfill\else
\agt, Volume \thevolumenumber\ (\thevolumeyear)\hfill\fi}
\def\@evenfoot{\@oddfoot}
\makeatother
\fi
\let\maketitlepage\makeagttitle
\let\makeshorttitle\maketitlepage
\let\maketitle\maketitlepage

\newwrite\gtoutfile
\long\gdef\makeheadfile{  
{\def\\{, }\def\s{ }
\immediate\openout\gtoutfile head.xxx
\immediate\write\gtoutfile{To: math@arxiv.org}
\immediate\write\gtoutfile{Subject: put OR rep NNNNN:ppppp}
\immediate\write\gtoutfile{--text follows this line--}
\immediate\write\gtoutfile{Proxy-for: \ifx\theasciiauthors\relax
\theauthors\else\theasciiauthors\fi\s<\ifx\theasciiemail\relax\theemail\else\theasciiemail\fi>}
\immediate\write\gtoutfile{\noexpand\\}
\immediate\write\gtoutfile{Authors: \ifx\theasciiauthors\relax
\theauthors\else\theasciiauthors\fi}
{\def\\{ }\immediate\write\gtoutfile{Title: \ifx\theasciititle\relax
\thetitle\else\theasciititle\fi}}
\immediate\write\gtoutfile{Subj-class: GT or SG, GR etc}
\immediate\write\gtoutfile{MSC-class: \theprimaryclass\ifx\thesecondaryclass\relax\else, \thesecondaryclass\fi}
\immediate\write\gtoutfile{Journal-ref: Algebr. Geom. Topol. \thevolumenumber\s
(\thevolumeyear) \startpage-\finishpage}
\immediate\write\gtoutfile{Comments: Published by Algebraic and
Geometric Topology at}
\immediate\write\gtoutfile{\s\s\s  http://www.maths.warwick.ac.uk/agt/AGTVol\thevolumenumber/agt-\thevolumenumber-\thepapernumber.abs.html}
\immediate\write\gtoutfile{\noexpand\\}
\immediate\write\gtoutfile{}
\ifx\theasciiabstract\relax
\immediate\write\gtoutfile{\theabstract}\else
\immediate\write\gtoutfile{\theasciiabstract}\fi
\immediate\write\gtoutfile{}
\immediate\write\gtoutfile{\noexpand\\}
\immediate\write\gtoutfile{}
\immediate\closeout\gtoutfile}}  

\def\maketitlepage{\makeagttitle\makeheadfile}
\let\makeshorttitle\maketitlepage
\let\maketitle\maketitlepage

\def\ifplaintex{\expandafter\ifx\csname documentclass\endcsname\relax}

\def\gtp{{\mathsurround=0pt\it $\cal G\mskip-2mu$eometry \&\ 
$\cal T\!\!$opology $\cal P\!$ublications}}  

\def\recd{{\small Received:\qua\receiveddate\ifx\reviseddate\relax
\else\qquad Revised:\qua\reviseddate\fi\par}} 

\def\lognumber#1{\def\thelognumber{#1}}
\def\volumenumber#1{\def\thevolumenumber{#1}}
\def\volumeyear#1{\def\thevolumeyear{#1}}
\def\papernumber#1{\def\thepapernumber{#1}}
\def\pagenumbers#1#2{\def\startpage{#1}\def\finishpage{#2}}
\def\published#1{\def\publishdate{#1}}

\def\received#1{\def\receiveddate{#1}}

\def\accepted#1{\def\accepteddate{#1}}
\def\asciititle#1{\def\theasciititle{#1}}
\def\covertitle#1{\def\thecovertitle{#1}}

\def\asciiaddress#1{\def\theasciiaddress{#1}}

\long\def\asciiabstract#1{\long\def\theasciiabstract{#1}}
\def\asciikeywords#1{\def\theasciikeywords{#1}}

\let\\\par\let\thelognumber\relax\let\thevolumenumber\relax
\let\thepapernumber\relax\let\thevolumeyear\relax\let\startpage\relax
\let\finishpage\relax\let\publishdate\relax\let\receiveddate\relax
\let\reviseddate\relax\let\accepteddate\relax\let\theasciititle\relax
\let\thecovertitle\relax\let\theasciiauthors\relax\let\theasciiaddress\relax
\let\theasciiabstract\relax\let\theasciikeywords\relax

\let\theasciiemail\relax

\ifplaintex
\font\logobig=cmssbx10 scaled 3836
\font\logomed=cmssbx10 scaled 2557
\else
\font\logobig=cmssbx10 scaled 4200
\font\logomed=cmssbx10 scaled 2800
\fi

\long\def\makeagttitle{   
\count0=\startpage
\agt\hfill      
\hbox to 45truept{\vbox to 0pt{\vglue -13truept{\logomed A\kern -.37em{\logobig 
T}\kern -.38em G}\vss}\hss}
\break
{\small Volume \thevolumenumber\ (\thevolumeyear)
\startpage--\finishpage\nl
Published: \publishdate}

\vglue .25truein

{\parskip=0pt\leftskip 0pt plus
1fil\def\\{\par\smallskip}{\Large\bf\thetitle}\par\medskip} \vglue
0.05truein

{\parskip=0pt\leftskip 0pt plus 1fil\def\\{\par}{\sc\theauthors}
\par\medskip}%
 
\vglue 0.03truein 

{\small\leftskip 25truept\rightskip 25truept{\bf Abstract}\stdspace\theabstract

{\bf AMS Classification}\stdspace\theprimaryclass
\ifx\thesecondaryclass\relax\else; \thesecondaryclass\fi\par
{\bf Keywords}\stdspace \thekeywords\par}\vglue 7truept

}   

\ifplaintex
\font\phead=cmsl9 scaled 950
\font\pnum=cmbx10 scaled 913
\font\pfoot=cmsl9 scaled 950
\headline{\vbox to 0pt{\vskip -4.5mm\line{\small\phead\ifnum
\count0=\startpage ISSN 1472-2739 (on-line) 1472-2747 (printed)
\hfill {\pnum\folio}\else\ifodd\count0\def\\{ }%
\ifx\theshorttitle\relax\thetitle\else\theshorttitle\fi\hfill{\pnum\folio}
\else\def\\{ and }{\pnum\folio}\hfill\ifx\theshortauthors\relax\theauthors
\else\theshortauthors\fi\fi\fi}\vss}}
\footline{\vbox to 0pt{\vglue 0mm\line{\small\pfoot\ifnum\count0=\startpage
\copyright\ \gtp\hfill\else
\agt, Volume \thevolumenumber\ (\thevolumeyear)\hfill\fi}\vss}}
\else
\font=cmsl9 scaled 1050
\font\lnum=cmbx10 
\font=cmsl9 scaled 1050
\makeatletter
\def\@oddhead{{\small\lhead\ifnum\count0=\startpage ISSN 1472-2739 
(on-line) 1472-2747 (printed)\hfill {\lnum\number\count0}\else\ifodd\count0
\def\\{ }\ifx\theshorttitle\relax \thetitle \else\theshorttitle\fi\hfill
{\lnum\number\count0}\else\def\\{ and }{\lnum\number\count0}
\hfill\ifx\theshortauthors\relax 
\theauthors\else\theshortauthors\fi\fi\fi}}\def\@evenhead{\@oddhead}
\def\@oddfoot{\small\lfoot\ifnum\count0=\startpage\copyright\ \gtp\hfill\else
\agt, Volume \thevolumenumber\ (\thevolumeyear)\hfill\fi}
\def\@evenfoot{\@oddfoot}
\makeatother
\fi
\let\maketitlepage\makeagttitle
\let\makeshorttitle\maketitlepage
\let\maketitle\maketitlepage

\newwrite\gtoutfile
\long\gdef\makeheadfile{  
{\def\\{, }\def\s{ }
\immediate\openout\gtoutfile head.xxx
\immediate\write\gtoutfile{To: math@arxiv.org}
\immediate\write\gtoutfile{Subject: put OR rep NNNNN:ppppp}
\immediate\write\gtoutfile{--text follows this line--}
\immediate\write\gtoutfile{Proxy-for: \ifx\theasciiauthors\relax
\theauthors\else\theasciiauthors\fi\s<\ifx\theasciiemail\relax\theemail\else\theasciiemail\fi>}
\immediate\write\gtoutfile{\noexpand\\}
\immediate\write\gtoutfile{Authors: \ifx\theasciiauthors\relax
\theauthors\else\theasciiauthors\fi}
{\def\\{ }\immediate\write\gtoutfile{Title: \ifx\theasciititle\relax
\thetitle\else\theasciititle\fi}}
\immediate\write\gtoutfile{Subj-class: GT or SG, GR etc}
\immediate\write\gtoutfile{MSC-class: \theprimaryclass\ifx\thesecondaryclass\relax\else, \thesecondaryclass\fi}
\immediate\write\gtoutfile{Journal-ref: Algebr. Geom. Topol. \thevolumenumber\s
(\thevolumeyear) \startpage-\finishpage}
\immediate\write\gtoutfile{Comments: Published by Algebraic and
Geometric Topology at}
\immediate\write\gtoutfile{\s\s\s  http://www.maths.warwick.ac.uk/agt/AGTVol\thevolumenumber/agt-\thevolumenumber-\thepapernumber.abs.html}
\immediate\write\gtoutfile{\noexpand\\}
\immediate\write\gtoutfile{}
\ifx\theasciiabstract\relax
\immediate\write\gtoutfile{\theabstract}\else
\immediate\write\gtoutfile{\theasciiabstract}\fi
\immediate\write\gtoutfile{}
\immediate\write\gtoutfile{\noexpand\\}
\immediate\write\gtoutfile{}
\immediate\closeout\gtoutfile}}  

\def\maketitlepage{\makeagttitle\makeheadfile}
\let\makeshorttitle\maketitlepage
\let\maketitle\maketitlepage

\lognumber{27}
\volumenumber{2}
\volumeyear{2002}
\papernumber{27}
\published{20 July 2002}
\pagenumbers{563}{590}
\received{24 January 2002}
\accepted{10 July 2002}

\usepackage{amsmath,amssymb}

\newtheorem{theorem}{Theorem}[section]
\newtheorem{corollary}[theorem]{Corollary}
\newtheorem{proposition}[theorem]{Proposition}
\newtheorem{lemma}[theorem]{Lemma}

\newtheorem{example}[theorem]{Example}
\newtheorem{remark}[theorem]{Remark}
\newtheorem{definition}[theorem]{Definition}

\def\Hq#1{\widehat{P}^{#1}}

\begin{document}
\title{Linking first occurrence polynomials\\over ${\mathbb F}_p$ by
Steenrod operations}
\covertitle{Linking first occurrence polynomials\\over ${\noexpand\bf F}_p$ by
Steenrod operations}
\asciititle{Linking first occurrence polynomials over F_p by
Steenrod operations}

\author{Pham Anh Minh\\Grant Walker}
\email{paminh@dng.vnn.vn, grant@ma.man.ac.uk}

\address{
Department of Mathematics, College of Sciences\\University of Hue,
Dai hoc Khoa hoc, Hue, Vietnam\\{\rm and}\\Department of Mathematics, 
University of Manchester\\Oxford Road,
Manchester M13 9PL, U.K. }

\asciiaddress{
 Department of Mathematics, College of Sciences\\University of Hue,
Dai hoc Khoa hoc, Hue, Vietnam\\and\\Department of Mathematics, 
University of Manchester\\Oxford Road,
Manchester M13 9PL, U.K. }

\begin{abstract} 
This paper provides analogues of the results of \cite{linking} for odd
primes $p$. It is proved that for certain irreducible representations
$L(\lambda)$ of the full matrix semigroup $M_n({\mathbb F}_p)$, the
first occurrence of $L(\lambda)$ as a composition factor in the
polynomial algebra ${\bf P}={\mathbb F}_p[x_1, \ldots, x_n]$ is linked
by a Steenrod operation to the first occurrence of $L(\lambda)$ as a
submodule in ${\bf P}$. This operation is given explicitly as the
image of an admissible monomial in the Steenrod algebra ${\cal A}_p$
under the canonical anti-automorphism $\chi$. The first occurrences of
both kinds are also linked to higher degree occurrences of
$L(\lambda)$ by elements of the Milnor basis of ${\cal A}_p$.
\end{abstract}

\asciiabstract{ 
This paper provides analogues of the results of [G.Walker and
R.M.W. Wood, Linking first occurrence polynomials over F_2 by Steenrod
operations, J. Algebra 246 (2001), 739--760] for odd primes p. It is
proved that for certain irreducible representations L(\lambda) of the
full matrix semigroup M_n(F_p), the first occurrence of L(\lambda) as
a composition factor in the polynomial algebra P=F_p[x_1,...,x_n] is
linked by a Steenrod operation to the first occurrence of L(\lambda)
as a submodule in P. This operation is given explicitly as the image
of an admissible monomial in the Steenrod algebra A_p under the
canonical anti-automorphism \chi. The first occurrences of both kinds
are also linked to higher degree occurrences of L(\lambda) by elements
of the Milnor basis of A_p.}

\asciikeywords{Steenrod algebra, anti-automorphism, p-truncated
polynomial algebra T, T-regular partition/representation} 
\maketitle

\section{Introduction} \label{intro}

Our aim  is to obtain results corresponding to those of \cite{linking} for the case where the prime $p
>2$.  In this we are only partly successful. The main theorem of \cite{linking} gives a Steenrod operation which
links the first occurrence of each irreducible representation $L(\lambda)$ of the full matrix semigroup
$M_n({\mathbb F}_2)$ in the polynomial algebra ${\bf P}={\mathbb F}_2[x_1,  \ldots, x_n]$ with the first
occurrence of $L(\lambda)$ as a submodule in ${\bf P}$.  Here  $M_n({\mathbb F}_2)$ acts on ${\bf P}$ on the right
by linear substitutions, which commute with the action of the Steenrod algebra ${\cal A}_2$ on  ${\bf P}$ on the
left. By `first occurrence' we have in mind the decomposition  ${\bf P}= \sum_{d \ge 0}{\bf P}^d$, where ${\bf
P}^d$ is the module of homogeneous polynomials of total degree $d$, and the known facts that there are minimum
degrees $d_c(\lambda)$ and $d_s(\lambda)$ in which $L(\lambda)$ occurs, uniquely in each case, as a composition
factor and as a submodule respectively.

For an odd prime $p$,  we have again the commuting actions of $M_n =M_n({\mathbb F}_p)$ on the right of the
polynomial algebra ${\bf P}={\mathbb F}_p[x_1,  \ldots, x_n ]$ and the algebra ${\cal A}_p$ of Steenrod $p$th
powers (no Bocksteins) on the left.  We refer to ${\cal A}_p$, somewhat inaccurately, as the Steenrod algebra, and
grade it so that $P^r$ raises degree by $r(p-1)$.  There are  $p^n$
isomorphism classes  of irreducible ${\mathbb F}_p[M_n]$-modules
$L(\lambda)$,  indexed by  partitions $\lambda = (\lambda_1,
\lambda_2, \ldots, \lambda_n)$, which are column $p$-regular, i.e.\
$0 \le \lambda_i - \lambda_{i+1} \le p-1$ for $1 \le i \le n$,
where $\lambda_{n+1}=0$ \cite{Green,Harris-Kuhn,James-Kerber}. The
problem solved in \cite{linking}  is certainly more difficult in
this context.  The submodule degree $d_s(\lambda)$ has recently
been determined \cite{Minh-Tri} for every  irreducible
representation $L(\lambda)$ of $M_n$, but $d_c(\lambda)$ is not
known in general. In particular, the first occurrence problem
appears to be difficult even for the 1-dimensional representations
$\det^k$, $1 \le k \le p-3$, $p>3$, see \cite{Carlisle-Kuhn,
Carlisle-Walker}, although it is solved for $\det^{p-2}$
\cite{Carlisle85}.   (The partition indexing $\det^k$ is $(k,
\ldots, k)= (k^n)$, i.e.\ $k$ repeated $n$ times.)  Further, it is
not known in general whether ${\bf P}^{d_c(\lambda)}$ has a unique
composition factor isomorphic to  $L(\lambda)$. Here we identify a
class of irreducible representations $L(\lambda)$ which behave
systematically. Since they arise  naturally by considering tensor
powers of the $p$-truncated polynomial algebra ${\bf T} = {\bf
P}/(x_1^p, \ldots, x_n^p)$, we call them  {\em ${\bf T}$-regular}.

Our main result, Theorem \ref{chim}, gives  a Steenrod operation
$\theta(\lambda)$ which links  the first occurrence and the first
submodule occurrence in ${\bf P}$ of a {\bf T}-regular
$L(\lambda)$. This determines $d_c(\lambda)$ in the {\bf
T}-regular case. The operation $\theta(\lambda)$ is given
explicitly as the image of an admissible monomial under the
canonical anti-automorphism $\chi$ of ${\cal A}_p$. Calculations
for $n \le 3$ suggest that such an operation $\theta(\lambda)$ may
exist for every irreducible representation $L(\lambda)$ of $M_n$,
but we do not pursue this here. Tri \cite{Tri00} has given an
`algebraic' alternative to this `topological' method of finding
$d_c(\lambda)$, using coefficient functions of ${\mathbb
F}_p[M_n]$-modules.

 For $p=2$, ${\bf T}$ may be identified with the exterior algebra
$\Lambda(x_1, \ldots, x_n)$, and  all the irreducible
representations $L(\lambda)$ of $M_n$ are ${\bf T}$-regular.  For
$p>2$, the only  irreducible 1-dimensional ${\bf T}$-regular
representations of $M_n$ are the `trivial' representation, in
which all matrices act as 1, and the $\det^{p-1}$ representation,
in which non-singular matrices act as 1 and singular matrices as
0.  The `trivial' representation, for which $\lambda = (0)$,
occurs in ${\bf P}$ only as ${\bf P}^0$, the constant polynomials.
Our key example is the $\det^{p-1}$ representation. This occurs
first as a composition factor as the top degree ${\bf T}^{n(p-1)}$
of ${\bf T}$, where it is generated  by the monomial $(x_1x_2
\cdots x_n)^{p-1}$ modulo $p$th powers, and first as a submodule
in degree $p_n =(p^n-1)/(p-1)$, where it is generated by the
Vandermonde determinant
 $$
  w(n) = \left| \begin{array}{llll} x_1 & x_2 & \cdots & x_n \\
x_1^p & x_2^p & \cdots & x_n^p \\ \vdots & \vdots & \ddots &
\vdots \\ x_1^{p^{n-1}} & x_2^{p^{n-1}} & \cdots & x_n^{p^{n-1}}
\end{array} \right|.
$$
\begin{theorem} \label{detp-1}
Let $\chi$ be the  canonical anti-isomorphism of ${\cal
A}_p$. Then for\break    $n \ge 1$,
 $$
  \chi(P^{p_n-n})(x_1x_2 \cdots x_n)^{p-1} =  w(n)^{p-1},
 $$
 where $p_n=(p^n-1)/(p-1)$.
\end{theorem}

This result is true for $p=2$ if we interpret $P^r$ as $Sq^r$
\cite{linking}. The operation $\chi(P^{p_n-n})$  may be replaced
by the admissible monomial $P^{p^{n-1}-1}\cdots P^{p^2-1}P^{p-1}$,
which is identical to  the Milnor basis element $P(p-1,  \ldots,
p-1)$ of length $n-1$ (see Proposition \ref{hatrel}). In general
the operation $\chi(P^{r_1}P^{r_2}\cdots P^{r_m})$ used in Theorem
\ref{chim} can not be replaced by an admissible monomial or a
Milnor basis element.

The structure of the paper is as follows.  Section \ref{prelim}
contains basic facts about the action of $\chi(P^r)$ and Milnor
basis elements on polynomials.   Section \ref{basiccase} contains
independent proofs of Theorem \ref{detp-1} using invariant theory
and  by  direct computation. In Section \ref{background} we
introduce the class of ${\bf T}$-regular partitions  to which our
main results apply, and extend Theorem \ref{detp-1} to ${\bf T}^d$
for all $d$. The main results are stated in Section \ref{results}
and proved in Section \ref{proofs}.  Section \ref{reps}  relates
these results to the  ${\mathbb F}_p[M_n]$-module structure of
${\bf P}$.   Section \ref{Milnorstuff} gives Milnor basis elements
which link  the first occurrence and (in certain cases) the first
submodule occurrence of a ${\bf T}$-regular representation of
$M_n$ with submodules in higher degrees.

The  remarks which follow are intended to place our results in
topological, combinatorial and algebraic contexts.  As for
topology,    recall (e.g.\ \cite{Barcelona}) that there is an
${\cal A}_p$-module decomposition ${\bf P} = \oplus_\lambda
\delta(\lambda){\bf P}(\lambda)$, where the $\lambda$-isotypical
summand ${\bf P}(\lambda)$ is an indecomposable ${\cal
A}_p$-module, and where $\delta(\lambda)= \dim L(\lambda)$, the
dimension of $L(\lambda)$.  Identifying  ${\bf P}$ with the
cohomology algebra $H^*({\mathbb C}P^\infty \times \cdots \times
{\mathbb C}P^\infty; {\mathbb F}_p)$,  this decomposition can be
realized (after localization at  $p$) by a homotopy equivalence
$\Sigma({\mathbb C}P^\infty \times \cdots \times {\mathbb
C}P^\infty) \sim \bigvee_\lambda \delta(\lambda) Y_\lambda$, which
splits the suspension of the product of $n$ copies of infinite
complex projective space ${\mathbb C}P^\infty$ as a topological
sum of spaces $Y_\lambda$ such that $H^*(Y_\lambda; {\mathbb F}_p)
= \Sigma{\bf P}(\lambda)$. The family of ${\cal A}_p$-modules ${\bf
P}(\lambda)$ is of major interest in algebraic topology. From this
point of view, we determine the connectivity of $Y_\lambda$ for
{\bf T}-regular $\lambda$ (Corollary \ref{firstocc}) and find a
nonzero cohomology operation $\theta(\lambda)$ on its bottom class
(Theorem \ref{chim}).

 As for combinatorics and algebra,  our aim is to provide
information relating the ${\cal A}_p$-module structure of ${\bf
P}(\lambda)$ to combinatorial properties of $\lambda$ and
representation theoretic properties of  $L(\lambda)$. The
operation $\theta(\lambda)$ and its source and target polynomials
are combinatorially determined by $\lambda$. The target polynomial
is defined by $w(\lambda') = \prod_{j=1}^{\lambda_1}
w(\lambda'_j)$, where $\lambda'$ is the conjugate of $\lambda$, so
that $w(\lambda')$ is a product of determinants corresponding to
the columns of the diagram of $\lambda$. This polynomial has
already appeared in various forms in the literature. In Green's
description \cite[(5.4d)]{Green} of the highest weight vector of
the dual Weyl module $H^0(\lambda)$, $w(\lambda')$ appears as a
`bideterminant' in the coordinate ring of $M_n(K)$, where $K$ is
an infinite field of characteristic $p$.  A proof that
$w(\lambda')$ generates a submodule of ${\bf P}^{d_s(\lambda)}$
isomorphic to $L(\lambda)$ was given in \cite[Proposition
1.3]{Doty-Walker96}, and a proof that this is the first occurrence
of $L(\lambda)$ as a submodule in ${\bf P}$ was given  in
\cite{Minh-Tri}.

We would like to thank the referee of this paper for a very
careful reading and for a number of helpful suggestions.  

\section{Preliminary results} \label{prelim}

In this section we use variants of the Cartan formula $ P^r (f g)
 = \sum_{r=s+t} P^s f \cdot  P^t g$ to study the action on
 polynomials of the elements $\chi(P^r)$  and Milnor basis
 elements $P(R)$ in the Steenrod algebra ${\cal A}_p$. We begin
 with the standard formula
\begin{equation} \label{evalP}
P^i(x^{p^b})=
 \begin{cases}
x^{p^{b+1}} & \text{if $i=p^b$},\\
 0 & \text{otherwise for $i>0$}.
	\end{cases}
\end{equation}
In particular, we wish to evaluate Steenrod operations on
Vandermonde determinants of the form
 $$
  [x_{i_1}^{s_1},x_{i_2}^{s_2}, \ldots ,x_{i_n}^{s_n}] = \left|
\begin{array}{cccc}
x_{i_1}^{s_1} & x_{i_2}^{s_1} & \ldots & x_{i_n}^{s_1} \\
x_{i_1}^{s_2} & x_{i_2}^{s_2} & \ldots & x_{i_n}^{s_2} \\ \vdots &
\vdots & \ddots &\vdots \\ x_{i_1}^{s_n} & x_{i_2}^{s_n} & \ldots
& x_{i_n}^{s_n} \\
\end{array} \right|,
$$
 where the exponents $s_1, \ldots, s_n$ are powers of $p$.   As
above, we  shall abbreviate such determinants  by listing their
diagonal entries in square brackets: in particular, $w(n) = [x_1,
x_2^p, \ldots, x_n^{p^{n-1}}]$.  As in Theorem \ref{detp-1},  we
write $p_n = (p^n-1)/(p-1)$, so that $p_0=0$ and $p_n-p_j =
(p^n-p^j)/(p-1)$.   The following result is a straightforward
calculation using the Cartan formula and (\ref{evalP}).
\begin{lemma}  \label{Minhlemma}
If $r =  p_n-p_  j$, $0 \le j \le n$, then $$P^r w(n) = [x_1,
x_2^p, \ldots, x_j^{p^{j-1}}, x_{j+1}^{p^{j+1}}, \ldots,
x_n^{p^n}], $$ and $P^r w(n) =0$  otherwise. In particular,  $P^r
w(n) =0$  for $0 < r < p^{n-1}$. \qed
\end{lemma}

 To simplify signs, we usually write $\Hq{r}$ for
$(-1)^r\chi(P^r)$.  Thus if $v$ is one of the generators $x_i$ of
${\bf P}$, or more generally any linear form $v=\sum_{i=1}^n
a_ix_i$ in ${\bf P}^1$,
\begin{equation}  \label{hatx}
\Hq{r}v = 
\begin{cases}
v^{p^b} & \text{if $r=p_b$, $b \ge 0$},\\
 0 & \text{otherwise}.
\end{cases}
\end{equation}

  Formula (\ref{hatx}) follows from (\ref{evalP})  by using the
identity $\sum_{i+j=r} (-1)^i P^i \Hq{j} =0$ in ${\cal A}_p$ and
induction on $r$. Using the identity $\sum_{i+j=r} (-1)^i \Hq{i}
P^j =0$ and induction on $k$, (\ref{hatx}) can be generalized to
\begin{equation}  \label{hatx2k}
\Hq{r}x^{p^k} =
\begin{cases}
x^{p^b}  & \text{if $r= p_b-p_k$, $b \ge k$}, \\
0 & \text{otherwise}.
\end{cases}
\end{equation}
This leads to the following generalization of  \cite[Lemma
2.2]{linking}.
\begin{lemma} \label{vdM1}
$$ \Hq{r}[x_1^{p^k}, x_2^{p^{k+1}}, \ldots, x_n^{p^{k+n-1}}] =
\begin{cases}
[x_1^{p^k},  \ldots, x_{n-1}^{p^{k+n-2}},x_n^{p^b} ] & \text{if
$r=p_b-p_{k+n-1}$}, \\
 0 & \text{otherwise}.      \end{cases} $$
\end{lemma}
The modifications required to the proof given in \cite{linking}
are straightforward.  \endproof

In evaluating the operations $\Hq{r}$, we shall frequently make
use of the Cartan formula expansion for polynomials $f, g \in {\bf
P}$:
\begin{equation} \label{coprod}
\Hq{r}(fg) = \sum_{s+t=r} \Hq{s}f \cdot \Hq{t}g,
\end{equation}
 which holds because  $\chi$  is a  coalgebra homomorphism.

\begin{lemma} \label{coprod2} For all polynomials $f, g$ in ${\bf P}$ and all $r \ge 0$,
 $$
  \Hq{r}(f^pg) = \sum_{r=ps+t} (\Hq{s}f)^p \Hq{t}g.
  $$
\end{lemma}
\proof By (\ref{coprod}) it suffices to prove the case
$g=1$, i.e.\
 $$
  \Hq{r}f^p = 
\begin{cases}
 (\Hq{s}f)^p  & \text{if $r=ps$},\\
 0 & \text{if $r$ is    not divisible by $p$}.
\end{cases}
 $$
  In this case, the Cartan formula (\ref{coprod}) gives $\Hq{r}f^p
   = \sum \Hq{r_1}f \cdots \Hq{r_p}f$, where the sum is over all
   ordered decompositions $r=\sum_{i=1}^p r_i$, $r_i \ge 0$.
   Except in the case where $r_1= \ldots = r_p=s$, cyclic
   permutation of $r_1, \ldots,r_p$ gives $p$ equal terms which
   cancel in the sum.  \endproof

We write $\alpha(k)$ for the sum of the digits in the base $p$
expansion of a positive integer $k$, i.e.\ if $k = \sum_{i \ge 0}
a_i p^i$ where $0 \le a_i \le p-1$, then $\alpha(k)= \sum_{i \ge
0} a_i$.  Thus $\alpha(k)$ is the minimum number of powers of $p$
which have sum $k$, and $\alpha(k) \equiv k$ mod $p-1$. Formula
(\ref{hatx}) leads to the following simple sufficient condition
for the vanishing of $\Hq{r}$ on a homogeneous polynomial of
degree $d$.

\begin{lemma} \label{alphatest}
If   $\alpha(r(p-1)+d)>d$, then $\Hq{r}f=0$ for all  $f \in {\bf
P}^d$.
\end{lemma}

{\bf Proof}\qua Since the action of $\Hq{r}$ is linear and
commutes with specialization of the variables, it is sufficient to
prove this when  $f= x_1x_2 \cdots x_d$.  By  (\ref{coprod}) $
\Hq{r}f = \sum\Hq{r_1}x_1\Hq{r_2}x_2 \cdots \Hq{r_d}x_d$, where
the sum is over all ordered decompositions $r=r_1+r_2+\ldots + r_d
$ with $r_1, r_2, \ldots, r_d \ge 0$. By (\ref{hatx}), the only
non-zero terms are those in which $r_i = p_{k_i}$ for some
non-negative integers $k_1, k_2, \ldots, k_d$. But then $r(p-1)+d
= \sum_i p^{k_i}$, and the result follows by definition of
$\alpha$.  \endproof

\begin{lemma} \label{x^{p-1}}
Let  $k \ge 0$ and let  $v=\sum_{i=1}^n a_ix_i$ be a linear form
in ${\bf P}^1$. Then
$$
 \Hq{p^k-1}v^{p-1} =v^{p^k(p-1)}.
  $$
\end{lemma}

\proof There is a unique way to write $p^k-1$ as the
sum of $p-1$ integers of the form $p_i$ for $i \ge 0$, namely
$p^k-1 = (p-1)p_k$.  The result now follows from (\ref{hatx}) and
the Cartan formula (\ref{coprod}).  \endproof

\begin{remark} \label{x^{p-1}remark} {\rm
The same method can be used to evaluate $\Hq{r}v^{p-1}$ for all
$r$. The result is $$ \Hq{r}v^{p-1} = 
\begin{cases}
 c_r v^{(r+1)(p-1)} & \text{if $\alpha((r+1)(p-1)) = p-1$}, \\
 0 & \text{otherwise},
\end{cases}
 $$ where if
$(r+1)(p-1) = j_1 p^{a_1}+ \ldots + j_s p^{a_s}$, with  $a_1 >
\ldots > a_s \ge 0$ and $\sum_{i=1}^s j_i = p-1$,  then
$c_r=(p-1)!/(j_1!j_2!\cdots j_s!)$. }
\end{remark}

The following result, the  `Cartan formula for Milnor basis
elements' is well-known (cf. \cite[Lemma 5.3]{linking}).
\begin{lemma} \label{Cartan} For a Milnor basis
element $P(R)= P(r_1, \ldots, r_n)$ and polynomials $f, g \in {\bf
P}$,
 $$
  P(R) (f g) = \sum_{R=S+T} P(S)f \cdot P(T)g,
 $$
  where the sum is over all sequences $S= (s_1, \ldots, s_n)$ and
$T=(t_1, \ldots, t_n)$ of non-negative integers such that $r_i =
s_i+t_i$ for $1 \le i \le n$.  \qed
\end{lemma}
In the same way as for Lemma \ref{coprod2},  this gives  the
following result.
\begin{lemma} \label{Cartan3} Let  $P(R)=P(r_1, \ldots, r_n)$ be
a Milnor basis element and let $f, g \in {\bf P}$ be
polynomials. Then
$$
 P(R) (f^p g) = \sum_{R=pS+T} (P(S)f)^p \cdot P(T)g. \eqno{\qed}
 $$
\end{lemma}
 Here $R=pS+T$  means that $r_i = ps_i+t_i$ for $1 \le i \le n$.

\section{The $\det^{p-1}$ representation} \label{basiccase}

  In this section we give three proofs of Theorem \ref{detp-1}.
The first uses the results of \cite{Minh-Tri} on submodules,
while the second is a variant of this which uses only classical
invariant theory. The third proof is computational. The  first two
proofs use  the following preliminary result, which shows that the
operation $\Hq{p_n-n}$ maps to $0$ all  monomials  of  degree
$n(p-1)$ other than the generating monomial  $(x_1x_2 \cdots
x_n)^{p-1}$ for $\det^{p-1}$.

\begin{lemma}  \label{divbyx^p}
Let $f$ be a polynomial in  ${\bf P}^{n(p-1)}$ which is divisible
by $x^p$ for some variable $x=x_i$, $1 \le i \le n$. Then
$\Hq{p_n-n}f=0$.
\end{lemma}
{\bf Proof}\qua Let $f=x^p g$, where $g \in {\bf P}$. Then by
Lemma \ref{coprod2}
\begin{equation} \label{r(n)}
\Hq{p_n-n}f = \sum_{p_n-n=ps+t} (\Hq{s}x)^p \Hq{t}g.
\end{equation}
By  (\ref{hatx}), $\Hq{s}x =0$ if $s \neq p_k$ for some $k$ with
$0 \le k \le n-2$. Thus it is sufficient to prove that $\Hq{t}g
=0$ for $t= p_n-n-p\cdot p_k$, where $g \in {\bf P}^{n(p-1)-p}$.
By Lemma \ref{alphatest}, this holds when $\alpha((t+n)(p-1)-p) >
n(p-1)-p$. Now $(t+ n)(p-1)-p = p_n(p-1) -p\cdot
p_k(p-1)-p=p^n-p^{k+1}-1$, hence $\alpha((t+ n)(p-1)-p) = n(p-1)-1
> n(p-1)-p$ as required.  Thus $\Hq{t}g=0$ in all    terms of
(\ref{r(n)}) in which $\Hq{s}x \neq 0$, and so  $\Hq{p_n-n}f=0$.
\endproof

{\bf  First Proof of Theorem \ref{detp-1}}\qua We  first show
that the monomial   $m=(x_1 x_2^p \cdots$ $x_n^{p^{n-1}})^{p-1}$
appears in $\Hq{p_n-n} (x_1 \cdots x_n)^{p-1}$ with coefficient
1. In the Cartan formula expansion (\ref{coprod}), $m$ can appear
only in the term arising from the decomposition $p_n-n = r_1+r_2+
\ldots + r_n$, where $r_k = p^{k-1}-1$ for $1 \le k \le n$. By
Lemma \ref{x^{p-1}},  $m$ appears in this term with coefficient 1.

By  Lemma \ref{divbyx^p}, $\Hq{p_n-n}$ maps all other monomials in
degree $n(p-1)$ to $0$. Hence   $\Hq{p_n-n} (x_1 \cdots
x_n)^{p-1}$ generates a 1-dimensional ${\mathbb
F}_p[M_n]$-submodule of ${\bf P}^{p^n-1}$.  Since $(x_1 \cdots
x_n)^{p-1}$ generates the 1-dimensional quotient ${\bf
T}^{n(p-1)}$ of ${\bf P}^{n(p-1)}$ and  since ${\bf T}^{n(p-1)} \cong
\det^{p-1}$, this submodule of ${\bf P}^{p^n-1}$ is also
isomorphic to $\det^{p-1}$.

 It is known \cite{Minh-Tri} that the first submodule occurrence
of  $\det^{p-1}$  for $M_n$ in ${\bf P}$ is generated by
$w(n)^{p-1}$, and that this is the unique submodule occurrence of
 $\det^{p-1}$ in degree $p^n-1$. Since $m$ is the
product of the leading diagonal terms in $w(n)^{p-1} = [x_1,
x_2^p, \ldots, x_n^{p^{n-1}}]^{p-1}$, $m$ also has coefficient 1
in $w(n)^{p-1}$.  \endproof

{\bf  Second Proof of Theorem \ref{detp-1}}\qua We recall that
 $D(n,p)$ is the ring of $GL_n({\mathbb F}_p)$-invariants in {\bf
 P}, and that it is a polynomial algebra over ${\mathbb F}_p$ with
 generators $Q_{n,i}$ in degree $p^n-p^i$ for $0 \le i \le
 n-1$. We may identify  $Q_{n,0}$ with $w(n)^{p-1}$. Since  ${\bf
 T}^{n(p-1)}$ is isomorphic to the trivial $GL_n({\mathbb
 F}_p)$-module, it follows as in our first proof that $\Hq{p_n-n}
 (x_1 \cdots x_n)^{p-1} \in D(n,p)$.

We shall prove that $w(n)$ divides $\Hq{p_n-n} (x_1 \cdots
x_n)^{p-1}$.  Recall that $w(n)$ is the product of linear factors
$c_1 x_1 + \ldots + c_n x_n$, where $c_1, \ldots, c_n \in {\mathbb
F}_p$. If we specialize the variables in  $(x_1 \cdots x_n)^{p-1}$
by imposing the relation $c_1 x_1 + \ldots + c_n x_n =0$, then
every monomial in the resulting polynomial is divisible by $x^p$
for some variable $x=x_i$. By Lemma \ref{divbyx^p}, such a
monomial is in the kernel of $\Hq{p_n-n}$. Thus $\Hq{p_n-n} (x_1
\cdots x_n)^{p-1}$ is divisible by $c_1 x_1 + \ldots + c_n x_n$,
and so it is divisible by $w(n)$.

Now an  element of $D(n,p)$ in degree $p^n-1$ which is divisible
by $w(n)$ must be a scalar multiple of $Q_{n,0} = w(n)^{p-1}$.
For if a polynomial in the remaining generators $Q_{n,1}, \ldots,
Q_{n,n-1}$ of $D(n,p)$ is divisible by $w(n)$, the quotient would
be  $SL_n({\mathbb F}_p)$-invariant,
 giving a non-trivial polynomial relation between
$Q_{n,1}, \ldots, Q_{n,n-1}$ and $w(n)$.  This contradicts
Dickson's theorem that these  are algebraically independent
generators of the polynomial algebra of $SL_n({\mathbb
F}_p)$-invariants in {\bf P}. \endproof

Our third proof of Theorem  \ref{detp-1} is by direct calculation.
We shall evaluate  the Milnor basis element $P(p-1, \ldots, p-1)$
of length $n-1$ on $(x_1 \cdots x_n)^{p-1}$.  The following result
relates   the element $P(p-1, \ldots, p-1,b)$ of length $n$ to
admissible monomials and to the anti-automorphism $\chi$. In
particular, we show that $P(p-1, \ldots, p-1)$ and $\Hq{p_n-n}$
have the same action on  $(x_1 \cdots x_n)^{p-1}$.
\begin{proposition} \label{hatrel} For  $1 \le b \le
p-1$,
\begin{itemize}
\item[\rm(i)] $P(p-1, \ldots, p-1,b) =  P^{(b+1)p^{n-1}-1} \cdots
P^{(b+1)p-1}P^b$ for  $n \ge 1$,
\item[\rm(ii)]  $\Hq{(b+1)p_n-n}g = P(p-1, \ldots, p-1,b)g$ if $\deg
g \le  n(p-1)+b$ for $n \ge 1$,
\item[\rm(iii)] $ \Hq{(b+1)p_n-n} = P^{(b+1)p^{n-1}}
\Hq{(b+1)p_{n-1}-n} +P(p-1, \ldots, p-1,b)$ for  $n \ge 2$. 
\end{itemize}
\end{proposition}
\proof
Statement {(i)} is a special case of
 \cite[Theorem 1.1]{Monks}.  For   {(ii)}, recall
 \cite{Milnor} that $\Hq{d}$ is the sum of all Milnor basis
 elements $P(R)$ in  degree $d(p-1)$.  Here  $R=(r_1, r_2,
 \ldots)$ is a finite sequence of non-negative integers, and
 $P(R)$ has degree $|R|= \sum(p^i-1)r_i$ and excess $e(R)= \sum
 r_i$.  In particular, $P^d = P(d)$ is the unique Milnor basis
 element of maximum excess $d$ in degree $d(p-1)$, but in  general
 there may be more than one  element of minimum excess in a given
 degree.

 We will show that $P(p-1, \ldots, p-1,b)$ is the unique  element
  of minimum excess $e=(n-1)(p-1)+b$ in degree $d(p-1)$ when $d=
  (b+1)p_n-n$.  By \cite[Lemma 8]{Milnor}  a bijection $P(r_1,
  r_2, \ldots, r_m) \leftrightarrow P^{t_1}P^{t_2} \cdots P^{t_m}$
  between the Milnor basis and the admissible basis of ${\cal
  A}_p$ is defined by $t_m=r_m$ and $t_i =r_i +pt_{i+1}$ for $1
  \le i <m$.  This  preserves both the degree and the excess. Thus
  it is equivalent to prove that $m= P^{(b+1)p^{n-1}-1} \cdots
  P^{(b+1)p-1}P^b$ is the unique admissible monomial of minimum
  excess  in degree $d(p-1)$. Now the excess of an admissible
  monomial $P^{t_1}P^{t_2} \cdots P^{t_m}$ is $pt_1 -d(p-1)$ where
  $d= \sum_i t_i$, and so it is minimal when $t_1$ is minimal.
  It is easy to verify that $m$ is the unique admissible monomial
  in degree $d(p-1)$  for which $t_1 = (b+1)p^{n-1}-1$, and that
  this value of $t_1$ is minimal.

 Note that $p$ divides $|R|+e(R)$ for all $R$.  Hence
$\Hq{(b+1)p_n-n}-P(p-1, \ldots, p-1, b)$ has excess $> e+p-1 =
n(p-1)+b$, and so $\Hq{(b+1)p_n-n}g = P(p-1, \ldots, p-1, b)g$
when $g$ is a polynomial of degree $\le n(p-1)+b$.

\medskip
 {(iii)}\qua Recall  Davis's formula \cite{Davis}
\begin{equation} \label{uvhat}
P^u\Hq{v} = \sum_{|R|=(p-1)(u+v)} \binom {|R|+e(R)}{pu} P(R),
\end{equation}
which we may apply   in the case $u= (b+1)p^{n-1}$, $v= (b+1)p_{n-1}-n$ to
show that $P^u\Hq{v}$ is the sum of all Milnor basis elements in
degree $d(p-1)$ other than   the  element  $P(p-1, \ldots,
p-1,b)$ of minimal excess.

For $R = (p-1, \ldots, p-1,b)$ we have  $|R|+e(R)= (b+1)p^n -p$, and
since $pu= (b+1)p^n$ the coefficient in (\ref{uvhat}) is zero. Since
$p$ divides $|R|+e(R)$ for all $R$,  $|R|+e(R) \ge (b+1)p^n$ for all
other $R$ with $|R|=d(p-1)$. As remarked above,  the unique
element of maximal excess is  $P^d$ itself, and so for all $R$ we
have $|R|+e(R) \le pd = (b+1)(p+p^2+ \ldots + p^n) -pn$. It is clear
from this inequality that the coefficient in (\ref{uvhat}) is 1
for all $R \neq (p-1, \ldots, p-1,b)$.  \endproof

{\bf Third Proof of Theorem \ref{detp-1}}\qua Let $\theta_{n} =
P^{p^n-1} \cdots P^{p^2-1}P^{p-1}$ for $n \ge 1$, and  $\theta_0
=1$. We assume  that $\theta_{n-1}(x_1 \cdots x_n)^{p-1} =
w(n)^{p-1}$ as induction hypothesis on $n$, the case $n=1$ being
trivial.

 The cofactor expansion of  $w(n+1) = [x_1, x_2^p, \ldots,
  x_{n+1}^{p^n}]$ by the top row gives $w(n+1) = \sum_{i=1}^{n+1}
  (-1)^i x_i \Delta_i^p$, where $\Delta_i = [x_1,  \ldots,
  x_{i-1}^{p^{i-2}}, x_{i+1}^{p^{i-1}}, \ldots,
  x_{n+1}^{p^{n-1}}]$.  Hence $w(n+1)\cdot  (x_1 \cdots
  x_{n+1})^{p-1} = \sum_{i=1}^{n+1} (-1)^i  x_i^p \Delta_i^p(x_1
  \cdots x_{i-1} x_{i+1} \cdots x_{n+1})^{p-1}$.

  By Proposition \ref{hatrel}{(i)},  $\theta_n = P(p-1,
  \ldots, p-1)$ of length $n$, and so by Lemma \ref{Cartan3}
  $\theta_n(w(n+1)\cdot  (x_1 \cdots x_{n+1})^{p-1}) =
  \sum_{i=1}^{n+1}  (-1)^i x_i^p \Delta_i^p \theta_n (x_1 \cdots
  x_{i-1}      x_{i+1} \cdots x_{n+1})^{p-1}$. Since $\theta_n =
  P^{p^n-1}\theta_{n-1}$, $\theta_n (x_1 \cdots  x_{i-1}
  x_{i+1}\cdots x_{n+1})^{p-1} = P^{p^n-1}\Delta_i^{p-1}$ by the
  induction hypothesis.  Since $\Delta_i^{p-1}$ has degree
  $p^n-1$, $P^{p^n-1} \Delta_i^{p-1} = \Delta_i^{p(p-1)}$. Hence
  $\theta_n (w(n+1)\cdot  (x_1 \cdots x_{n+1})^{p-1}) =
  \sum_{i=1}^{n+1} (-1)^i x_i^p \Delta_i^{p^2} = w(n+1)^p$.

By Lemma \ref{Minhlemma},  $P^r w(n+1) = 0$ for $0 < r < p^{n}$.
As $\theta_n = P^{p^n-1} \cdots P^{p^2-1}P^{p-1}$, iterated
application of the Cartan formula gives $  \theta_n (w(n+1)\cdot
(x_1 \cdots x_{n+1})^{p-1}) = w(n+1)\cdot  \theta_n(x_1 \cdots
x_{n+1})^{p-1}$.  Hence $w(n+1)\cdot  \theta_n(x_1 \cdots
x_{n+1})^{p-1} = w(n+1)^p$. Cancelling the factor $w(n+1)$, the
inductive step is proved.

\section{${\bf T}$-regular partitions} \label{background}

In this section we  define the special class of {\em ${\bf
T}$-regular} partitions, and extend Theorem \ref{detp-1} to give a
Steenrod operation $\Hq{r}$ which links the first occurrence and
first submodule occurrence of  ${\bf T}^d$ for all $d$.  In fact
we prove a more general result which links the first occurrence to
a family of higher degree occurrences.

The truncated polynomial module ${\bf T}^d = {\bf P}^d/({\bf
P}^d\cap (x_1^p, \ldots, x_n^p))$ has a ${\mathbb F}_p$-basis
represented in ${\bf P}^d$ by the set of all monomials
$x_1^{s_1}x_2^{s_2} \cdots x_n^{s_n}$ of total degree $d = \sum_i
s_i$ with $s_i <p$ for $1 \le i \le n$.  By \cite[Theorem
6.1]{Carlisle-Kuhn} ${\bf T}^d \cong L((p-1)^{n-1} b)$, where
$d=(n-1)(p-1)+b$ and $1 \le b \le p-1$.  We regard the
corresponding diagram as a block of $p-1$ columns, in which the
first $b$ columns have length $n$ and the remaining  $p-b-1$
columns have length $n-1$.  Given a partition $\lambda$, we can
divide its diagram into $m$ blocks of $p-1$ columns and compare
the blocks with the diagrams corresponding to these.  (The $m$th
block may have $< p-1$ columns.)  For $1 \le j \le m$,  let
$\lambda_{(j)}$ be the partition whose diagram is the $j$th block,
and let $\gamma_j = \deg \lambda_{(j)}$ be the  number of boxes in
the $j$th block.

\begin{definition} \label{Treg} {\rm
 A  column $p$-regular partition $\lambda$ is {\em {\bf
T}-regular} if $L(\lambda_{(j)}) \cong {\bf T}^{\gamma_j}$ for all $j$.
Equivalently, for all $a \ge 1$, there is at most one value of $i$
for which $(a-1)(p-1) < \lambda_i < a(p-1)$. If $\lambda$ is {\bf
T}-regular, we call $\gamma$ the {\em {\bf T}-conjugate} of
$\lambda$. }
\end{definition}

 In the case $p=2$, all column $2$-regular partitions are ${\bf
T}$-regular, and $\gamma =\lambda'$, the conjugate of $\lambda$.
If $\kappa$ is column 2-regular, then  the partition
$\lambda=(p-1)\kappa$ obtained by multiplying each part of
$\kappa$ by $p-1$ is ${\bf T}$-regular. Since  $\lambda$ is column
$p$-regular, $\gamma_j -\gamma_{j+1} \ge p-1$ for all $j$, and  $m
\le n$.  Thus there is a bijection $\lambda \leftrightarrow
\gamma$ between the set of  {\bf T}-regular partitions $\lambda=
(\lambda_1, \ldots, \lambda_n)$  and the set of partitions
$\gamma= (\gamma_1, \ldots, \gamma_n)$ which satisfy $\gamma_1\le
n(p-1) $ and $\gamma_j -\gamma_{j+1} \ge p-1$ for $1 \le j \le
n-1$. In terms of  the Mullineux involution $M$ on the set of all
row $p$-regular partitions, $\lambda$ and $\gamma$ are related by
$M(\gamma)= \lambda'$ \cite[Proposition 3.13]{MSF}.

  We next  extend Theorem \ref{detp-1} to give linking formulae
   for  the representations ${\bf T}^d$.  It will be convenient to
   introduce abbreviated notation for some further Vandermonde
   determinants.  Let $w(n,a) = [x_1, \ldots, x_a^{p^{a-1}},
   x_{a+1}^{p^{a+1}}, \ldots, x_n^{p^n}]$  for $0 \le a \le n$,
   where the exponent $p^a$ is omitted. In particular, $w(n,n) =
   w(n)$ and $w(n,0) = w(n)^p$.
\begin{proposition} \label{Minhtrick2}
For  $n \ge 1$ and   $1 \le i \le p-1$, let $i= i_1+ \cdots +i_s$
where $i_1, \ldots, i_s >0$, and let $j= i_1p_{a_1}+ \cdots +
i_sp_{a_s}$, where $a_1> \ldots  > a_s \ge 0$. Then
 $$
 \Hq{p_n-n-j}\left((x_1x_2 \cdots x_{n-1})^{p-1}
  x_n^{p-i-1}\right) = (-1)^{i(n-1)-j}w(n)^{p-i-1} \cdot
  \prod_{r=1}^s w(n-1,a_r)^{i_r}.
 $$
\end{proposition}

Specializing to the case $s=1$, $j=ip_{n-1}$ and putting $b=
p-1-i$, we obtain an operation  linking the first occurrence and
the first submodule occurrence of  the representation ${\bf T}^d$,
as follows. Theorem \ref{detp-1} can be taken as the case $b=0$ or as the case
$b=p-1$; we choose $b=p-1$ to fit notation later. 

\begin{corollary} \label{Minhtrick}
For  $n \ge 1$ and  $1 \le b \le p-1$,
 $$
  \Hq{(b+1)p_{n-1}-(n-1)}\left((x_1x_2 \cdots x_{n-1})^{p-1}
x_n^b\right) = w(n)^b \cdot w(n-1)^{p-b-1}.
 $$
\end{corollary}

{\bf  Proof of Proposition \ref{Minhtrick2}}\qua We introduce a
  parameter into Theorem \ref{detp-1}, by working in  ${\mathbb
  F}_p[x_1,\ldots, x_{n+1}]$ and writing $x_{n+1}=t$ in order to
  distinguish this variable. Since the action of ${\cal A}_p$
  commutes with the linear substitution which maps $x_n$ to
  $x_n+t$ and fixes $x_i$ for $i \neq n$,   we obtain
\begin{equation} \label{x_n+t}
\Hq{p_n-n}(x_1 \cdots x_{n-1}(x_n+t))^{p-1} = [x_1, x_2^p, \ldots,
x_{n-1}^{p^{n-2}}, (x_n+t)^{p^{n-1}}]^{p-1}.
\end{equation}
Expanding the left hand side of (\ref{x_n+t}) by the binomial
theorem, we obtain
 $$
  \sum_{i=0}^{p-1} (-1)^i \Hq{p_n-n}((x_1 \cdots x_{n-1})^{p-1}
  x_n^{p-1-i} t^i).
 $$
  The right hand side of (\ref{x_n+t}) is
   $$
    [x_1, x_2^p, \ldots, x_{n-1}^{p^{n-2}},
x_n^{p^{n-1}} \hspace{-5pt}+t^{p^{n-1}}]^{p-1} \hspace{-3pt}
 = \sum_{i=0}^{p-1} (-1)^i
w(n)^{p-1-i}  [x_1, x_2^p, \ldots, x_{n-1}^{p^{n-2}},
t^{p^{n-1}}]^i,
 $$
  since $w(n) = [x_1, x_2^p, \ldots, x_n^{p^{n-1}}]$.  The
   summands in (\ref{x_n+t}) corresponding to $i=0$  give the
   original   result, Theorem \ref{detp-1}, and so are equal.  In
   fact  we can  equate the $i$th summands for all $i$.  This
   happens because  $\Hq{r}$ raises degree by $r(p-1)$, so that
   the powers  $t^k$ which occur in the $i$th summand on the left
   have $k \equiv i$ mod $p-1$, while if $t^k$ occurs in the $i$th
   summand on the right, then  $k$ is the sum of  $i$ powers of
   $p$, so that again $k \equiv i$ mod $p-1$.  Hence for $1 \le i
   \le p-1$ we have \begin{equation}  \label{icomp}
   \Hq{p_n-n}((x_1 \cdots x_{n-1})^{p-1} x_n^{p-1-i} t^i) =
   w(n)^{p-1-i} \cdot [x_1, x_2^p, \ldots, x_{n-1}^{p^{n-2}},
   t^{p^{n-1}}]^i.  \end{equation} Since the powers $t^k$ of $t$
   which can appear here  are such that $k$ is the sum of $i$
   powers of $p$, we can write $k = i_1 p^{a_1} + \ldots + i_s
   p^{a_s}$, where $a_1> \ldots > a_s \ge 0$  and $i_1 + \ldots +
   i_s = i$.  Using the expansion
 $$
  [x_1, x_2^p, \ldots, x_{n-1}^{p^{n-2}}, t^{p^{n-1}}] =
  \sum_{a=0}^{n-1} (-1)^{n-1-a} w(n-1,a) t^{p^a}
$$
 we can evaluate the coefficient of $t^k$ on the right hand side
 of (\ref{icomp}) as
 $$
   (-1)^{i(n-1)-j} \frac{i!}{i_1! \cdots i_s!} w(n)^{p-1-i} \cdot
   w(n-1,a_1)^{i_1} \cdots w(n-s,a_s)^{i_s},
 $$
where we have simplified the sign by noting that $a_1i_1+ \ldots
   +a_si_s \equiv j$ mod $2$  since $p_a \equiv a$ mod 2,.  By the
   Cartan formula (\ref{coprod}), the left hand side of
   (\ref{icomp}) is
 $$
  \sum_{j=0}^{p_n-n} \Hq{p_n-n-j}\left((x_1 \cdots x_{n-1})^{p-1}
  x_n^{p-1-i}\right) \cdot \Hq{j}t^i
 $$
Here the term in $t^k$ arises from $\Hq{j}t^i$ where $k=j(p-1)+i$,
so that $j = i_1 p_{a_1}+ \ldots +  i_s p_{a_s}$, and since this
decomposition of $j$ as a sum of at most $i$ powers of $p$ is
unique, formulas (\ref{hatx}) and (\ref{coprod}) give $  \Hq{j}
t^i =(i!/ i_1! \cdots i_s!) t^k$.  Thus equating coefficients of
$t^k$ in  (\ref{icomp}) gives the result.

\section{Linking  for {\bf T}-regular representations}
\label{results}

In this section we state our main results. We fix an odd prime $p$
and a positive integer $n$ throughout. As in \cite{linking},  our
results will be statements about polynomials in $n$ variables when
$\lambda$ has length $n$, i.e.\  $\lambda$ has $n$ nonzero parts.
There is no loss of generality, since the projection in $M_n$
which sends $x_n$ to 0 and $x_i$ to $x_i$ for $i<n$ maps
$L(\lambda)$ to zero if $\lambda_n >0$ and on to the corresponding
${\mathbb  F}_p[M_{n-1}]$-module $L(\lambda)$ if $\lambda_n =0$
(cf. \cite[Section 3]{Carlisle-Kuhn}). Hence  we shall always
assume  that $\lambda_n \neq 0$.

  We first establish some notation. Given a ${\bf T}$-regular
partition $\lambda$ of length $n$, we define a polynomial
$v(\lambda)$ whose degree $d_c(\lambda)$  is given by  (\ref{dc})
and which `represents'  $L(\lambda)$, in the sense that the
submodule  of ${\bf P}^{d_c(\lambda)}$ generated by $v(\lambda)$
has  a quotient module isomorphic to $L(\lambda)$. We index the
diagram of $\lambda$  using  matrix  coordinates $(i,j)$,   so
that   $1 \le i \le n$ and $1 \le j \le \lambda_i$.
\begin{definition} \label{antidiag}
 The $k$th {\em antidiagonal} of the  diagram of $\lambda$ is
the set of boxes such that $j+i(p-1)=k+p-1$. If  the lowest box
is in row $i$ and the highest is in row $i-s+1$, let $
v_k(\lambda)= [x_{i-s+1}, x_{i-s+2}^p, \ldots, x_i^{p^{s-1}}]$,
and let $v(\lambda) = \prod_{k=1}^{\gamma_1}v_k(\lambda)$.
\end{definition}
 Thus an antidiagonal is the set of boxes which lie on a line of
  slope $1/(p-1)$ in the  diagram, and $v(\lambda)$ is a product
  of corresponding  Vandermonde determinants.  Indenting
  successive rows by $p-1$ columns, we obtain a shifted diagram
  whose  columns correspond to these antidiagonals. The ${\bf
  T}$-conjugate $\gamma$ of $\lambda$ records the number of
  antidiagonals $\gamma_s$ of length $\ge s$ for all $s \ge 1$.

\begin{example} \label{antidiagex}
{\rm Let $p=5$, $\lambda =(9,6,3)$, so that $\gamma =
(11,6,1)$. The shifted diagram $$\begin{array}{ccccccccccc}
*&*&*&*&*&*&*&*&*\\&&&&*&*&*&*&*&*\\  &&&&&&&&*&*&* \\
\end{array}$$
gives $v(\lambda) = x_1^4 \cdot [x_1,x_2^5]^4 \cdot
[x_1,x_2^5,x_3^{25}] \cdot[x_2,x_3^5] \cdot x_3$. }
\end{example}

Recall \cite{Minh-Tri} that $w(\lambda') = \prod_{j=1}^{\lambda_1}
w(\lambda'_j)$  generates the first occurrence of $L(\lambda)$ as
a submodule in ${\bf P}$. Thus we can rewrite the linking theorem
for  ${\bf T}^d$, Corollary \ref{Minhtrick}, as follows.

\begin{theorem} \label{Minhtrick3}  Let $d =(n-1)(p-1)+b$, where $n \ge 1$ and $1 \le b \le p-1$, so that ${\bf
T}^d \cong L(\lambda)$ where $\lambda = ((p-1)^{n-1}b)$.  Then  $
  \Hq{r}v(\lambda) = w(\lambda')$, where $r = (b+1)p_{n-1}-(n-1)$
  and $p_{n-1}=(p^{n-1}-1)/(p-1)$. \qed \end{theorem}

By the  {\em leading monomial}  of a polynomial we mean  the
monomial  $\prod_{i=1}^n x_i^{s_i}$ occurring in it (ignoring the
nonzero coefficient) whose exponents are highest in left
lexicographic order.  The leading monomial $s(\lambda)$ of
$v(\lambda)$  is obtained by multiplying the principal
antidiagonals in the determinants $v_k(\lambda)$, $1 \le k \le
\gamma_1$.   (In Example \ref{antidiagex},
$s(\lambda)=x_1^{49}x_2^{14}x_3^3$.) The base $p$ expansion of
every exponent in $s(\lambda)$ has the form $s_i = c_k p^k +
(p-1)p^{k-1} + \ldots + (p-1)p + (p-1)$, i.e.\ $s_i \equiv
-1$~mod~$p^k$, where $p^k < s_i < p^{k+1}$. We adapt the
terminology introduced by Singer \cite{Singer}, by calling such a
monomial a `spike'. In the case $p=2$, $s(\lambda) =
x_1^{2^{\lambda_1}-1}\cdots x_n^{2^{\lambda_n}-1}$. A polynomial
which contains such a spike can not be `hit', i.e.\ it can not be
the image of a polynomial of lower degree under a Steenrod
operation. This is easily seen by considering the 1-variable
case. Hence the polynomial $v(\lambda)$ is not hit.

\begin{proposition} \label{spike}   Let $\lambda$ be ${\bf T}$-regular with ${\bf T}$-conjugate $\gamma$.
 \begin{itemize}
 \item[\rm(i)] If $\lambda_i = a_i(p-1)+b_i$,
  $a_i \ge 0$, $1 \le b_i \le p-1$, then $ s(\lambda) =
  \prod_{i=1}^n x_i^{(b_i+1)p^{a_i}-1}$.
\item[\rm(ii)] With  $\lambda_{(j)}$  as in Definition \ref{Treg}, $
s(\lambda) = v(\lambda_{(1)}) \cdot v(\lambda_{(2)})^p \cdots
v(\lambda_{(m)})^{p^{m-1}}$.
  \item[\rm(iii)] The coefficient of
$s(\lambda)$ in $v(\lambda)$ is $(-1)^{\epsilon(\lambda)}$, with
$\epsilon(\lambda) = \sum_{j=1}^{[m/2]} (-1)^{j-1} \gamma_{2j}$.
\end{itemize}
\end{proposition}

{\bf Proof}\qua Formulae  {(i)} and {(ii)} are easily
read off from a tableau obtained by entering $p^{j-1}$ in each box
in the $j$th block of $p-1$ columns of the diagram of $\lambda$,
and reading this according to rows and to blocks of columns.  For
{(iii)}, note that the sign of the term arising from the
leading antidiagonal in the expansion of an  $s \times s$
determinant is $+1$ for $s \equiv 0,1$ mod $4$ and $-1$ for $s
\equiv 2,3$ mod $4$, and that the diagram of $\lambda$ has
$\gamma_j$  antidiagonals of length $\ge j$.  \endproof

In Theorem \ref{chir} we establish {(i)} a `level 0 formula',
which gives a sufficient condition for $\Hq{r}v(\lambda) =0$, and
{(ii)} a `level 1 formula', which  gives a sufficient
condition for $\Hq{r}v(\lambda)$ to be a product related to the
decomposition $\lambda = \lambda_{(1)} + \lambda^-$ which splits
off the first $p-1$ columns of the diagram. Thus  $\lambda_{(1)} =
((p-1)^{n-1}b)$, where $\gamma_1 =(n-1)(p-1)+b$ and $1 \le b \le
p-1$, and $\lambda^-$ is defined by $\lambda^-_i =
\lambda_i-(p-1)$ if $\lambda_i \ge p-1$, and $\lambda^-_i  =0$
otherwise.  Our main linking result, Theorem \ref{chim}, follows
from Theorem \ref{chir}  by induction on $m$, the length  of
$\gamma$.  The proofs  of Theorems \ref{chir} and \ref{chim} are
deferred to Section \ref{proofs}.

\begin{theorem} \label{chir}   Let $\lambda$ be    {\bf T}-regular
 with ${\bf T}$-conjugate $\gamma$, let $d_c$ be defined by {\rm
  (\ref{dc})} below, and let $R(r, \lambda)
  =r(p-1)+d_c(\lambda)-d_c(\lambda^-)$.  Recall that $\alpha(k)$
  is the sum of the digits in the base $p$ expansion of $k$.
\begin{itemize}
\item[\rm(i)] If  $\alpha(R(r, \lambda)) > \gamma_1$, then $\Hq{r}
v(\lambda) = 0$.
\item[\rm(ii)] If  $\alpha(R(r, \lambda)) = \gamma_1$, then $\Hq{r}
v(\lambda) = \Hq{r+d_c(\lambda^-)} v(\lambda_{(1)}) \cdot
v(\lambda^-)$.
\end{itemize}

\end{theorem}

\begin{remark} \label{chir2} {\rm
  Taking  $p=2$ and $P^r =Sq^r$, this reduces to \cite[Theorem
2.1]{linking}, since that theorem can be applied to
$\lambda_{(1)}=(1^n)$ to obtain   $\widehat{Sq}^{r+d_c(\lambda^-)}
v(\lambda_{(1)}) =[x_1^{2^{a_1}},  \ldots, x_n^{2^{a_n}}]$, where
$a_1 <  \ldots < a_n$. The hypothesis on $r$ is satisfied since
$r+d_c(\lambda^-)+n = r + d_c(\lambda)- d_c(\lambda^-) = 2^{a_1} +
\ldots + 2^{a_n}$. }
\end{remark}

Combining Theorem \ref{Minhtrick3}  with Theorem \ref{chir}, we
obtain our main theorem.

\begin{theorem} \label{chim} Let  $\lambda$ be  {\bf T}-regular
  with {\bf T}-conjugate  $\gamma$ of length $m$. For $1 \le k \le
m$, let $\gamma_k = (n_k-1)(p-1)+b_k$, where $n_k \ge 1$ and $1
\le b_k \le p-1$.  Then
 $$
 \Hq{r_m}\cdots \Hq{r_2}\Hq{r_1} v(\lambda) = w(\lambda'),
 $$
 where $r_k = (b_k+1) p_{n_k-1}-(n_k-1)-\sum_{j=k+1}^m p^{j-k-1}
 \gamma_j$.
\end{theorem}

This theorem determines the first occurrence degree
$d_c(\lambda)$ when $\lambda$ is ${\bf T}$-regular.

\begin{corollary} \label{firstocc} Let $\lambda$ be {\bf T}-regular with ${\bf T}$-conjugate $\gamma$.
 Then the degree  in which the irreducible module $L(\lambda)$
first occurs as a composition factor in the polynomial algebra
${\bf P}$ is given by \begin{equation} \label{dc} d_c(\lambda) =
\sum_{i=1}^m p^{i-1} \gamma_i, \end{equation} and the  ${\mathbb
F}_p[M_n]$-submodule of ${\bf P}^{d_c(\lambda)}$ generated by
$v(\lambda)$ has a quotient module isomorphic to $L(\lambda)$.
\end{corollary}
{\bf Proof}\qua By \cite{Doty-Walker96} or \cite{Minh-Tri}
$w(\lambda')$ generates a submodule of ${\bf P}^{d_s(\lambda)}$
isomorphic to $L(\lambda)$. By Theorem \ref{chim}, there is a
Steenrod operation $\theta=\theta(\lambda)$ and a polynomial
$v(\lambda) \in {\bf P}^d$, where $d$ is given by
(\ref{dc}), such that $\theta(v(\lambda)) = w(\lambda')$. Hence
the quotient of the submodule generated by $v(\lambda)$ in ${\bf
P}^d$ by the  intersection of this submodule with the kernel of
$\theta$ is a composition factor of ${\bf P}^d$ which is
isomorphic to $L(\lambda)$. Hence the first occurrence degree
$d_c(\lambda) \le d$. But $d_c(\lambda) \ge d$   by
\cite[Proposition 2.13]{Carlisle-Walker}, and hence $d_c(\lambda)
= d$.  \endproof

As an example,  for $p=3$ the partition $\lambda= (5,3,2)$ is {\bf
T}-regular with {\bf T}-conjugate $\gamma=(6,3,1)$. The  module
$L(5,3,2)$ first occurs as a composition factor in degree
$6+3\cdot 3+1\cdot 9 =24$, and  as a submodule in degree $5+3\cdot
3 +2\cdot 9 = 32$.  The calculations  of \cite{Carlisle85} and
\cite{Doty-Walker92} for $n \le 3$ support the conjecture that the
the first occurrence degree $d_c(\lambda)$ is given by the formula
above if and only if  $\lambda$ is ${\bf T}$-regular.

The  integers $r_i$ in Theorem \ref{chim} can be calculated from a
tableau ${\rm Tab}(\lambda)$ obtained by entering integers into
the diagram of $\lambda$ as follows: if a box in row $i$
is the highest box in its antidiagonal, write $p_{i-1}$ in that
box and continue down the antidiagonal, multiplying the number
entered at each step by $p$.

\begin{lemma} \label{rk} The sum of the numbers entered in the $k$th block of $p-1$ columns using the above rule is  $r_k$.
The element  $P^{r_1}P^{r_2}\cdots P^{r_m}$ is an admissible
monomial in  ${\cal A}_p$, i.e.\ $r_k \ge pr_{k+1}$ for $1 \le k
\le m -1$.
\end{lemma}

\begin{example} \label{65432} {\rm For $p=3$, $\lambda =
(6,5,4,3,2)$, we obtain $(r_1,r_2,r_3) = (100,20,1)$ using the tableau below.
  $$ {\rm Tab}(\lambda) =
\begin{array}{|c|c|c|c|c|c|}
\hline 0&0&0&0&0&0 \\ \hline 0&0&0&0&1 \\ \cline{1-5} 0&0&3&4 \\
\cline{1-4} 9&12&13 \\ \cline{1-3} 39&40  \\ \cline{1-2}
\end{array}
$$
 Noting that $\Hq{r} = (-1)^r\chi(P^r)$, in this case Theorem
 \ref{chim} states that in ${\bf P}^{300}$,
  $$
  \chi(P^{100} P^{20} P^1) \left(x_1^2 \cdot [x_1,x_2^3]^2 \cdot
[x_1,x_2^3,x_3^9]^2 \cdot [x_2,x_3^3,x_4^9] \cdot [x_3,x_4^3]
\cdot [x_4,x_5^3] \cdot x_5 \right)
$$
 $$
  = -[x_1,x_2^3,x_3^9,x_4^{27},x_5^{81}]^2\cdot
  [x_1,x_2^3,x_3^9,x_4^{27}]\cdot [x_1,x_2^3,x_3^9]\cdot
  [x_1,x_2^3]\cdot x_1.
  $$
   }
\end{example}

{\bf Proof of Lemma \ref{rk}}\qua The inequality  $r_k \ge
pr_{k+1}$ for $1 \le k \le m -1$ is clear from the algorithm, and
can also be checked directly from the definition of $r_k$. Since
$r_2(\lambda) = r_1(\lambda^-)$, and so on, we need only check
the algorithm  for $r_1$.

To do this, we introduce a second tableau by entering $p_{i-1}$ in
the $i$th row of the first block of $p-1$ columns and $-p^{j-2}$
in all the boxes in the $j$th block of $p-1$ columns for $j>1$. In
Example \ref{65432} this is as follows.
  $$
\begin{array}{|c|c|c|c|c|c|}
\hline 0&0&-1&-1&-3&-3 \\ \hline 1&1&-1&-1&-3 \\ \cline{1-5}
4&4&-1&-1 \\ \cline{1-4} 13&13&-1 \\ \cline{1-3} 40&40 \\
\cline{1-2}
\end{array}
$$
 The  entries in a antidiagonal running from the $(i,j)$ box for
$1 \le  j \le p-1$ are    then $p_{i-1}, -1, -p, \ldots,
-p^{s-2}$, and their sum  $p_{i-1}-p_{s-1} = p^{s-1}p_{i-s}$ is
the number entered in this box in ${\rm Tab}(\lambda)$.

It remains to check that the sum of all the entries in the second
tableau is $r_1= (b_1+1)p_{n-1}-(n-1)-d_c(\lambda^-)$. To see
this, note that the entries in $\lambda^-$ sum to
$-d_c(\lambda^-)$, while the entries in the last row of
$\lambda_{(1)}$ sum to $bp_{n-1}$ and the entries in the first
$n-1$ rows sum to $(p-1)(p_0+p_1+\ldots + p_{n-2}) =
p_{n-1}-(n-1)$.  \endproof

Since  $w(n)$ is a product of linear factors, so also is
$v(\lambda)$, and by   Theorems \ref{Minhtrick3} and
\ref{chir} so also is
$\Hq{r_1}v(\lambda)$.  The following calculation shows that
$v(\lambda)$ divides $\Hq{r_1}v(\lambda)$, and that the quotient
can be read off from ${\rm Tab}(\lambda)$  as follows:  replace
the entry $p_{i-1}-p_{s-1}$ in the $(i,j)$ box, $1 \le j \le p-1$,
by the product of all linear polynomials of the form $x_i +
\sum_{k<i} c_k x_k$, excluding those where $c_k=0$ for $1 \le k
\le i-s $.

\begin{corollary} \label{factors}
Let $\lambda$ be a {\bf T}-regular  partition. Let  the $k$th
antidiagonal in the diagram of $\lambda$ have length $s_k$ and
lowest box in row $n_k$.   Then
 $$
 \frac {\Hq{r_1} v(\lambda)}{ v(\lambda)} = \prod_{k=1}^{\gamma_1}
\prod_{\bf c} (c_1 x_1 + \ldots + c_{n_k-1}x_{n_k-1} + x_{n_k}),
 $$
where the inner product is over all vectors ${\bf c}= (c_1,
   \ldots, c_{n_k-1}) \in {\mathbb F}_p^{n_k-1}$ such that $(c_1,
   \ldots , c_{n_k-s_k}) \neq (0, \ldots ,0)$.
\end{corollary}

In Theorem \ref{detp-1}, $\lambda = ((p-1)^n)$, $v(\lambda) =
(x_1x_2 \cdots x_n)^{p-1}$ and $\Hq{r_1} v(\lambda) = [x_1, x_2^p,
\ldots, x_n^{p^{n-1}}]^{p-1}$. Since $s_k =1$ for $1 \le k \le
n(p-1)$, the quotient is the product of all linear polynomials in
$x_1,  \ldots, x_n$ which are not monomials.

\medskip
{\bf Proof of Corollary \ref{factors} }\qua The proof is by
induction on the number of antidiagonals $\gamma_1$.  Let $\phi(\lambda)
= \Hq{r_1} v(\lambda)/v(\lambda)$, where $r_1=r_1(\lambda)$. Let
$s$ denote the length of  the last
antidiagonal in the diagram of $\lambda$, and let $\mu$ be
the ${\bf T}$-regular partition obtained by removing this
antidiagonal  from the diagram of $\lambda$.  Then
by Theorems \ref{Minhtrick3} and \ref{chir},
 $$
\frac{\phi(\lambda)}{\phi(\mu)} = \frac{[x_1, x_2^p, \ldots,
x_n^{p^{n-1}}]}{ [x_1, x_2^p, \ldots, x_{n-1}^{p^{n-2}}]}
\cdot \frac{v(\lambda^-)}{ v(\mu^-)}   \cdot \frac{v(\mu)}{v(\lambda)}.
 $$
Note that $ \lambda^- = \mu^-$ when $s=1$.   Now  $[x_1,
x_2^p, \ldots, x_n^{p^{n-1}}]/ [x_1, 
x_2^p, \ldots, x_{n-1}^{p^{n-2}}] = \prod_{\bf c} (c_1 x_1 +
\ldots + c_{n-1}x_{n-1} + x_n)$, where the  product is taken over
all vectors ${\bf c}= (c_1, \ldots, c_{n-1}) \in {\mathbb
F}_p^{n-1}$.  Also $v(\lambda)/v(\mu) = v_{\gamma_1}(\lambda) =
[x_{n-s+1}, x_{n-s+2}^p, \ldots, x_n^{p^{s-1}}]$. Similarly
$v(\lambda^-)/v(\mu^-) = [x_{n-s+1}, x_{n-s+2}^p, \ldots,
x_{n-1}^{p^{s-2}}]$.  The quotient of these determinants is the
product of all $p^{s-1}$ linear polynomials  $ c_{n-s+1}
x_{n-s+1}+ \ldots + c_{n-1}x_{n-1}  +x_n$, so
$\phi(\lambda)/\phi(\mu) = \prod_{\bf c} (c_1 x_1 + \ldots +
c_{n-1}x_{n-1} +  x_n)$, where the  product is over all  ${\bf c}=
(c_1, \ldots, c_{n-1}) \in {\mathbb F}_p^{n-1}$ with $c_i \neq 0$
for some $i$ such that $1 \le i \le n-s$.

\section{Proof of the linking theorem} \label{proofs}

In this section we  prove Theorems \ref{chir} and \ref{chim}.  The
following  lemma will help in checking conditions on the numerical
function $\alpha$.
\begin{lemma} \label{alpha}
\begin{itemize}
\item[\rm(i)] Let $R \ge 1$ have base $p$ expansion  $R=j_1p^{a_1}+
\ldots + j_tp^{a_t}$, where $1 \le j_1, \ldots, j_t \le p-1$, $0
\le a_1 < \ldots < a_t$, and let  $k \ge 0$. Then $\alpha(R-p^k)
\ge \alpha(R)-1$, with equality if and only if $k=a_i$, $1 \le i
\le t$.
\item[\rm(ii)] With notation as in Theorem \ref{chir}, and with $\mu$
and $s$ as in the proof of  Corollary \ref{factors}, for $r \ge 1$
and $k \ge 0$ we have
 $$
R(r-p_k+p_{s-1},\mu) = R(r-p_k +d_c(\lambda^-), \mu_{(1)}) = R(r,
\lambda)-p^k.
 $$
\end{itemize}
\end{lemma}
{\bf Proof}\qua If $k \neq a_i$ for $1 \le i \le t$, then
subtraction of $p^k$  must yield at least one new term  $(p-1)p^a$
in the base $p$ expansion. This proves {\rm (i)}.   For {\rm
(ii)}, since $d_c(\lambda)=d_c(\lambda_{(1)})+pd_c(\lambda^-)$ and
$d_c(\lambda_{(1)})=\gamma_1$ we have $R= R(r, \lambda)=
(p-1)(r+d_c(\lambda^-))+\gamma_1$.  Comparing the first occurrence
degrees for $L(\lambda)$ and $L(\mu)$ given by (\ref{dc}),
\begin{equation} \label{comparedc}
 d_c(\lambda)= d_c(\mu) +p_s,
\quad d_c(\lambda^-) = d_c(\mu^-) +p_{s-1}, \quad
d_c(\lambda_{(1)})= d_c(\mu_{(1)}) +1. 
 \end{equation}
Hence we have
 $R(r-p_k+p_{s-1},\mu) = (p-1)(r-p_k+p_{s-1}+d_c(\mu^-)) +
d_c(\mu_{(1)}) = (p-1)(r-p_k+d_c(\lambda^-)) + d_c(\mu_{(1)})
=R(r-p_k +d_c(\lambda^-), \mu_{(1)}) = R-(p-1)p_k-1 =  R-p^k$.
\endproof

{\bf Proof of Theorem \ref{chir}(i)}\qua  We argue by induction
on $\gamma_1$, the number of antidiagonals  of $\lambda$.   With
$\mu$ and $s$ as above,  $v(\lambda) = [x_{n-s+1}, x_{n-s+2}^p,
\ldots, x_n^{p^{s-1}}]\cdot v(\mu)$.  Using formula (\ref{coprod})
and Lemma \ref{vdM1}, for all  $r \ge 1$ we have
\begin{equation} \label{split}
\Hq{r} v(\lambda) = \sum_{k \ge s-1} [x_{n-s+1}, x_{n-s+2}^p,
\ldots, x_{n-1}^{p^{s-2}}, x_n^{p^k}] \cdot \Hq{r-p_k +p_{s-1}}
v(\mu).
\end{equation}
By Lemma \ref{alpha}, if $\alpha(R(r, \lambda)) > \gamma_1$ then
$\alpha(R(r-p_k+p_{s-1},\mu)) > \gamma_1-1$ for all $k \ge 0$.
Since $\mu$ has $\gamma_1-1$ antidiagonals, the second factor in
each term of (\ref{split}) is zero by the induction hypothesis.
Hence $\Hq{r} v(\lambda) = 0$ if $\alpha(R(r,\lambda)) >\gamma_1$,
completing the induction.

\medskip
{\bf Proof of Theorem \ref{chir}(ii)}\qua As in Lemma
\ref{alpha}, let  $R= R(r,\lambda)$ have base $p$ expansion
$R=j_1p^{a_1}+ \ldots + j_tp^{a_t}$, let $\alpha(R)=\gamma_1$ and
let $R'= R(r-p_k+p_{s-1},\mu)$. Then the lemma gives $\alpha(R') =
\gamma_1-1$ if $k=a_i$, $1 \le i \le t$, and $\alpha(R') >
\gamma_1-1$ otherwise.  Hence,  applying part {\rm (i)} of the
theorem to  (\ref{split}), we have
$$
\Hq{r} v(\lambda) = \sum_{i=1}^t [x_{n-s+1}, x_{n-s+2}^p, \ldots,
x_{n-1}^{p^{s-2}}, x_n^{p^{a_i}}] \cdot \Hq{r-p_{a_i} +p_{s-1}}
v(\mu).
$$
 Since $\alpha(R(r-p_{a_i}+p_{s-1},\mu)=\gamma_1-1=d_c(\mu_{(1)})$
    by the lemma, and $p_{s-1}+d_c(\mu^-)=d_c(\lambda^-)$,    the
    inductive hypothesis on $\mu$ gives
 $$
  \Hq{r-p_{a_i}+p_{s-1}} v(\mu) =
  \Hq{r-p_{a_i}+d_c(\lambda^-)}v(\mu_{(1)}) \cdot v(\mu^-),\quad 1
  \le i \le t.
   $$
We can similarly  use the lemma  to simplify the right hand side
of the required identity. Since $v(\lambda_{(1)}) = x_n
v(\mu_{(1)})$, from (\ref{coprod}) and (\ref{hatx}) we have
 $$
   \Hq{r+d_c(\lambda^-)}v(\lambda_{(1)}) =\sum_{k \ge 0} x_n^{p^k}
   \Hq{r+d_c(\lambda^-)-p_k}v(\mu_{(1)}).
$$
 By the lemma, $R(r+d_c(\lambda^-)-p_k, \mu_{(1)}) = R-p^k$, so
 that by {\rm (i)} we can again reduce to the sum  over $k=a_i$,
 $1 \le i \le t$.  As $v(\lambda^-) =  [x_{n-s+1}, x_{n-s+2}^p,
 \ldots, x_{n-1}^{p^{s-2}}]\cdot  v(\mu^-)$,  it remains after
 cancelling the factor $v(\mu^-)$ and rearranging  terms to prove
 that
 $$
\sum_{i=1}^t \left([x_{n-s+1}, x_{n-s+2}^p, \ldots,
x_{n-1}^{p^{s-2}}, x_n^{p^{a_i}}] - [x_{n-s+1}, x_{n-s+2}^p,
\ldots, x_{n-1}^{p^{s-2}}] x_n^{p^{a_i}} \right) \cdot f_i =0,
 $$
  where $f_i= \Hq{r-p_{a_i}+d_c(\lambda^-)}v(\mu_{(1)})$.   The
expansion of the  $s \times s$ determinant in the $p^{a_i}$ powers
of the variables is
  $$
   \sum_{j=1}^s (-1)^{s-j} [x_{n-s+1}, \ldots,
x_{n-s+j-1}^{p^{j-2}},   x_{n-s+j+1}^{p^{j-1}}, \ldots,
x_n^{p^{s-2}}]x_{n-s+j}^{p^{a_i}}.
 $$
Thus the term with $j=s$ cancels, and interchanging the $i$ and
$j$ summations, the required formula becomes
 $$
 \sum_{j=1}^{s-1} (-1)^{s-j} [x_{n-s+1}, \ldots,
x_{n-s+j-1}^{p^{j-2}},   x_{n-s+j+1}^{p^{j-1}}, \ldots,
x_n^{p^{s-2}}] \cdot \sum_{i=1}^t x_{n-s+j}^{p^{a_i}}f_i =0.
 $$
Since $ \Hq{r+d_c(\lambda^-)}(x_{n-s+j}v(\mu_{(1)}))= \sum_{i=1}^t
 x_{n-s+j}^{p^{a_i}}f_i$ by a similar argument  using
 (\ref{coprod}), (\ref{evalP}) and Lemma \ref{alpha},  it suffices
 to prove that the monomial $x_{n-s+j} v(\mu_{(1)})$ is in the
 kernel of $\Hq{r+d_c(\lambda^-)}$ for $1 \le j \le s-1$.  This
 monomial  is divisible by $x_{n-s+j}^p$.  By permuting the
 variables, it suffices to consider the case where it is divisible
 by $x_1^p$.  Hence the proof of Theorem \ref{chir} is completed
 by the following calculation.

\begin{proposition} \label{zerocase}
Let $R=R(r,\lambda)$ and let $\alpha(R) = \gamma_1$, where
$\gamma_1 = (n-1)(p-1)+b$ and $ 1 \le b  \le p-1$. Then
 $$
  \Hq{r+d_c(\lambda^-)} (x_1^p ( x_2 \cdots x_{n-1})^{p-1} \cdot
  x_n^{b-1}) =0.
 $$
\end{proposition}

{\bf Proof}\qua By Lemma \ref{coprod2}, with $f=x_1$ and $g=  (
x_2 \cdots x_{n-1})^{p-1} \cdot  x_n^{b-1}$,
 $$
\Hq{u}(x_1^p \cdot g) = \sum_{u=pv+w} (\Hq{v}x_1)^p \cdot
\Hq{w}(g).
 $$
  Note that $g= v(\nu)$ where $\nu = ((p-1)^{n-2}(b-1))$. By
(\ref{hatx}), $\Hq{v}x_1=0$  for $v \neq p_k$, $k \ge 0$, so we
may assume that $w = u-pv = r+ d_c(\lambda^-) - p\cdot p_k$.
Since $p \cdot p_k = p_{k+1}-1$ and $d_c(\mu_{(1)})=
p-1+d_c(\nu)$, $R(w, \nu)= R(r- p_{k+1}+d_c(\lambda^-),
\mu_{(1)})=R -p^{k+1}$ by Lemma \ref{alpha}{\rm (ii)}.   Since
$\alpha(R) = \gamma_1$, Lemma \ref{alpha}{\rm (i)} gives
$\alpha(R(w, \nu)) \ge \gamma_1-1> \gamma_1-p$.   Since $d_c(\nu)
= \gamma_1-p$,  $\Hq{w}g=0$ by  Theorem \ref{chir}{\rm (i)}.
\endproof

{\bf Proof of Theorem \ref{chim}}\qua This follows from Theorem
  \ref{chir}  by induction on $m$.  Let $\gamma_1 = (n-1)(p-1)+b$,
  $1 \le b \le p-1$.  We wish to apply Theorem \ref{chir} with
  $r=r_1$, so we must check that $\alpha(R(r_1, \lambda)) =
  \gamma_1$.  For this, note that (\ref{dc}) gives $d_c(\lambda^-)
  = \sum_{j=2}^m p^{j-2}\gamma_j$, so that $r_1+d_c(\lambda^-)=
  (b+1)p_{n-1}-(n-1)$.  Thus $R(r_1, \lambda) =
  (p-1)(r_1+d_c(\lambda^-)) +\gamma_1 = (b+1)(p^{n-1}-1)
  -(p-1)(n-1)+\gamma_1 = b p^{n-1} + (p^{n-1}-1)$.  Hence $r_1$
  satisfies the hypothesis of Theorem \ref{chir}, so that
  $\Hq{r_1} v(\lambda) =  \Hq{r_1+d_c(\lambda^-)} v(\lambda_{(1)})
  \cdot v(\lambda^-)$.  By Theorem \ref{Minhtrick3},
  $\Hq{r_1+d_c(\lambda^-)} v(\lambda_{(1)}) = w(\lambda_{(1)}')$.

 Now $r_i(\lambda) = r_{i-1}(\lambda^-)$ for $2 \le i \le m$,  and
so the inductive step reduces to showing that \begin{equation}
\label{indstep} \Hq{r_m}\cdots \Hq{r_2}\left( w(\lambda'_{(1)})
\cdot v(\lambda^-)\right ) = w(\lambda'_{(1)}) \cdot
\Hq{r_m}\cdots \Hq{r_2} v(\lambda^-).  \end{equation} Recall from
Lemma \ref{rk} that  $r_1,  \ldots, r_m$ is an admissible
sequence, i.e.\ $r_k \ge pr_{k+1}$ for $k \ge 1$. Since $r_1 \le
(b+1)p_{n-1}$, $r_1 < p^{n-1}$ if $b < p-1$ and $r_1 < p^n$ if
$b=p-1$. Thus we can deduce (\ref{indstep}) from Lemma \ref{vdM1}
and the coproduct formula (\ref{coprod}), as follows. We have
$w(\lambda'_{(1)}) = w(n)^b w(n-1)^{p-1-b}$. Now $\Hq{r} w(n) =0$
for $0 < r < p^{n-1}$ and $\Hq{r} w(n-1) =0$ for $0 < r <
p^{n-2}$. If there are any factors $w(n-1)$ in
$w(\lambda'_{(1)})$, then $r_2 < p^{n-2}$, and otherwise it
suffices to have  $r_2 < p^{n-1}$.

\section{First occurrence submodules} \label{reps}

For a  ${\bf T}$-regular partition $\lambda$, the  ${\mathbb
  F}_p[M_n]$-submodule of ${\bf P}^{d_c(\lambda)}$ generated by
  the first occurrence polynomial $v(\lambda)$ is a
  `representative polynomial' for $L(\lambda)$ in the sense that
  this module has a quotient isomorphic to $L(\lambda)$(see
  Corollary \ref{firstocc}).  In the case where $\lambda =
  (p-1)\kappa$ for a column $2$-regular partition $\kappa$, the
  leading monomial  $s(\lambda)= x_1^{p^{\kappa_1}-1}\cdots
  x_n^{p^{\kappa_n}-1}$ has the same property.  This is implicit
  in the work of Carlisle and Kuhn \cite{Carlisle-Kuhn}, who
  identify a subquotient ${\bf T}^\gamma$ of ${\bf
  P}^{d_c(\lambda)}$ such that ${\bf T}^\gamma \cong {\bf
  T}^{\gamma_1}  \otimes \ldots \otimes {\bf T}^{\gamma_m}$, where
  $\gamma$ is the {\bf T}-conjugate of $\lambda$.  Explicitly, if
  $v_i \in {\bf T}^{\gamma_i}$ corresponds to a monomial in $x_1,
  \ldots, x_n$ with all exponents $< p$, then $v_1 \otimes \ldots
  \otimes v_m \in {\bf T}^{\gamma_1} \otimes \ldots \otimes {\bf
  T}^{\gamma_m}$ corresponds  to the equivalence class of $v_1
  \cdot v_2^p \cdots v_m^{p^{m-1}}$ in the appropriate subquotient
  of ${\bf P}^{d_c(\lambda)}$. Proposition \ref{spike}{\rm (ii)}
  shows that,  taking $v_j= v(\lambda_{(j)})$, this monomial is
  $s(\lambda)$. Tri \cite{Tri00} has recently proved that if
  $\lambda$ is {\bf T}-regular, then  $L(\lambda)$ is a
  composition factor in ${\bf T}^\gamma$.

We recall from \cite[Section 4]{linking} the notion  of a {\em
base $p$ $\omega$-vector}.
\begin{definition} \label{omegavector} {\rm
Given a prime   $p$,  the {\em base $p$  $\omega$-vector}
$\omega(s)$  of   a sequence
of non-negative integers $s=(s_1, \ldots, s_n)$   is
defined as follows. Write each $s_i$ in base $p$ as $s_i = \sum_{j\ge
1} s_{i,j}p^{j-1}$, where $0 \le s_{i,j} \le p-1$, and let $\omega_j(s) =
\sum_{i=1}^n  s_{i,j}$, i.e.\  add the base $p$ expansions
without `carries'. Then $\omega(s) = (\omega_1(s), \ldots, \omega_l(s))$,
with  {\em length} $l = \max \{j:\ \omega_j(s) > 0\}$ 
and {\em degree}  $d =
\sum_{i=1}^n s_i = \sum_{j= 1}^l \omega_j(s) p^{j-1}$.
}
\end{definition}

  Given $\omega$-vectors $\rho$ and
$\sigma$, we say that $\rho$ {\em
dominates} $\sigma$, and write $\rho \succeq \sigma$ or $\sigma
\preceq \rho$, if and only if $\sum_{i=1}^k p^{i-1} \rho_i \ge
\sum_{i=1}^k p^{i-1} \sigma_i$ for all $k \ge 1$. By the
$\omega$-vector of a monomial $\prod_{i=1}^n x_i^{s_i}$ we mean
the $\omega$-vector of its sequence of exponents  $s = (s_1,
\ldots, s_n)$. The dominance order on $\omega$-vectors of the same
degree is compatible with  left lexicographic order.

\begin{example} \label{omegalattice} {\rm
The lattice of base $p$ $\omega$-vectors of degree $1+p+p^2$ is
shown below. $$\begin{array}{ccccc} & &(1+p+p^2)&\\ && \downarrow
& \\ &&(1+p^2,1)&\\ && \downarrow & \\ &&\vdots &\\ & & \downarrow
& \\ &&(1+p,p)&\\ & \swarrow &&  \searrow \\ (1,1+p) &&&&
(1+p,0,1)\\ & \searrow && \swarrow \\ && (1,1,1) &
\end{array}$$}
\end{example}

\begin{proposition}  \label{Weylmod1} Let $\lambda$ be a {\bf T}-regular partition. Then the $\omega$-vector of the spike
monomial $s(\lambda)$ is the partition $\gamma$ {\bf T}-conjugate
to $\lambda$, and the polynomial $v(\lambda)$ is the sum of
$(-1)^{\epsilon(\lambda)}s(\lambda)$ and monomials $f$ such that
$\omega(f) \prec \gamma$.
\end{proposition}

{\bf Proof}\qua The proof is the same as that given in
\cite[proposition 4.5]{linking}, with $2$ replaced by $p$ and
$\lambda'$ replaced by $\gamma$.  For $\epsilon(\lambda)$, see
Proposition \ref{spike}{\rm (iii)}. \endproof

 Corollary \ref{firstocc} and Proposition \ref{Weylmod1} together
provide a `topological' proof that  the ${\mathbb
F}_p[M_n]$-submodule of ${\bf P}^{d_c(\lambda)}$ generated by
$s(\lambda)$ has a quotient module isomorphic to $L(\lambda)$. The
next result provides a further comparison between the spike
monomial $s(\lambda)$ and the polynomial $v(\lambda)$ in a special
case.  We conjecture that the corresponding statement holds for
all {\bf T}-regular partitions $\lambda$.

\begin{proposition} \label{bigsubmods} Assume that $\lambda_i = (p-1)\kappa_i$ for $1 \le i \le n$,
 where $\kappa = (\kappa_1, \ldots, \kappa_n)$ is a column
$2$-regular partition. Then the submodule of ${\bf
P}^{d_c(\lambda)}$ generated by the polynomial $v(\lambda)$ is
contained in the submodule generated by the  spike monomial
$s(\lambda)$.
\end{proposition}
The proof requires a preliminary lemma.
\begin{lemma} \label{LS}
If $f \in {\mathbb F}_p[x_2, \ldots, x_n]$ and $1 \le s \le n$,
then the ${\mathbb F}_p[M_n]$-submodule of ${\bf P}$ generated by
$x_1^{p^s-1}\cdot f$ contains $[x_1, x_2^p, \ldots,
x_s^{p^{s-1}}]^{p-1}\cdot f$.
\end{lemma}

{\bf Proof}\qua For each linear form $v = a_1 x_1 + \ldots + a_s
x_s$, where $a_i \in {\mathbb F}_p$ for $1 \le i \le s$, let $t_v:
{\bf P} \rightarrow {\bf P}$  be the transvection mapping $x_1$ to
$v$ and fixing $x_2, \ldots, x_n$.  We claim that the following
equation holds in ${\mathbb F}_p[x_1, \ldots, x_s]$.
\begin{equation} \label{sumI}
(-1)^s [x_1, x_2^p, \ldots, x_s^{p^{s-1}}]^{p-1} = \sum_v
v^{p^s-1}.
\end{equation}
 Since $t_v$ does not change the variables $x_2, \ldots, x_n$
which can occur in $f$, it follows from (\ref{sumI}) that $\sum_v
t_v$ is an element of the semigroup algebra ${\mathbb F}_p[M_n]$
which maps $x_1^{p^s-1}\cdot f$ to $(-1)^s[x_1, x_2^p, \ldots,
x_s^{p^{s-1}}]^{p-1}\cdot f$.

To prove (\ref{sumI}), first note that the  right hand side is
$GL_s({\mathbb F}_p)$-invariant. Further, it is mapped to $0$ by
every singular matrix $g \in M_s$, since  vectors   $(a_1, \ldots,
a_s)$ and  $(a'_1, \ldots, a'_s)$ in ${\mathbb F}_p^s$ in the same
coset of  the kernel of $g$ yield terms in (\ref{sumI}) with the
same image under $g$, and $p$ divides the order of this coset.
Arguing as in the first or second proof of Theorem \ref{detp-1},
with $s$ in place of $n$, it follows that (\ref{sumI}) holds up to
a (possibly zero) scalar.

Finally we verify that  the monomial $ m= x_1^{p-1}x_2^{p(p-1)}
\cdots x_s^{p^{s-1}(p-1)} $ has coefficient $(-1)^s$ in the right
hand side of (\ref{sumI}). For each linear form $v$, we have
$v^{p^s-1} = v^{p^{s-1}(p-1)}\cdots v^{p(p-1)} \cdot v^{p-1}$,
where  $ v^{p^j(p-1)}= (a_1 x_1^{p^j} + \ldots + a_s
x_s^{p^j})^{p-1}$ for $0 \le j \le s-1$. The exponent $p-1$ in $m$
must come from the last factor in this product, so we must choose
the term $(a_1x_1)^{p-1} = x_1^{p-1}$ from the last factor, and
$a_1 \neq 0$.  In the same way, we must choose the term
$(a_2x_2^p)^{p-1} = x_2^{p(p-1)}$ from the last but one factor,
and $a_2 \neq 0$.  Continuing in this way, we see that each of the
$(p-1)^s$ linear forms $v$ with all coefficients $a_i \neq 0$
gives a term containing $m$ (with coefficient 1), while other
choices of $v$ give terms not containing $m$. Thus the scalar
coefficient in (\ref{sumI}) is $(-1)^s$. \endproof

The following example shows how to  apply Lemma \ref{LS} to a
partition $\lambda$ of the form $(p-1)\kappa$, so as to generate
$v(\lambda)$ from $s(\lambda)$.
\begin{example} \label{spiketov} {\rm
Let $p=3$ and let  $\lambda = (6,6,4,4,2)$, so that  $s(\lambda) =
x^{26}y^{26}z^8t^8u^2$ and $v(\lambda) =
x^2[x,y^3]^2[x,y^3,z^9]^2[y,z^3,t^9]^2[t,u^3]^2$.

Begin by permuting the variables, so as to work with the spike
$u^8t^{26}z^{26}y^8x^2$.  Apply Lemma \ref{LS} with $x_1 =y$ and
$s=2$ to generate $ [y,x^3]^2 \cdot u^8t^{26}z^{26}x^2$. Repeat
with $x_1 = z$ and $s=3$ to generate $[z,y^3,x^9]^2\cdot
u^8t^{26}[y,x^3]^2x^2$,  then with $x_1=t$ and $s=3$ to generate
$[t,z^3,y^9]^2 \cdot u^8[z,y^3,x^9]^2[y,x^3]^2x^2$, and finally
with $x_1=u$ and $s=2$ to generate $v(\lambda)$. }
\end{example}

{\bf Proof of Proposition \ref{bigsubmods}}\qua We first observe
(see \cite[Proposition 4.9]{linking}) that the (multi)set of
lengths of the antidiagonals of the column $2$-regular partition
$\kappa$ is equal to the (multi)set of lengths of the rows.
Hence the spike monomial $ \tilde{s}(\lambda) = x_n^{p^{s_n}-1}
x_{n-1}^{p^{s_{n-1}-1}} \cdots x_1^{p^{s_1}-1}$, where $s_k$ is
the length of the $k$th antidiagonal of the diagram of $\kappa$,
can be obtained from $s(\lambda)$ by a suitable permutation of the
variables. We can now obtain $v(\lambda)$ from
$\tilde{s}(\lambda)$  by $n-1$ successive applications of Lemma
\ref{LS}, following the method illustrated by Example
\ref{spiketov}.

\section{${\bf T}$-regular partitions and the Milnor basis} \label{Milnorstuff}

In  this section we link the first occurrence polynomial
 $v(\lambda)$ and its leading monomial $s(\lambda)$ to the
 polynomial $p(\lambda') = \prod_{j=1}^m
 w(\lambda_{(j)}')^{p^{j-1}}$, which  generates a submodule
 occurrence of $L(\lambda)$ in a higher degree.   Here, as in
 Proposition \ref{spike}, $\lambda_{(j)}$ is   the partition given
 by the $j$th block of $p-1$ columns in the diagram of the ${\bf
 T}$-regular partition $\lambda$, and $m$ is the length of
 $\gamma$, the   ${\bf T}$-conjugate of $\lambda$.  In the case
 $\lambda = (p-1)\kappa$, we also link  the first submodule
 occurrence polynomial $w(\lambda')$ to $p(\lambda')$.  The
 linking is achieved by Milnor basis elements in ${\cal A}_p$
 which are combinatorially related to  $\lambda$. We also obtain a
 relation between monomials in ${\bf P}$ and Milnor basis elements
 in terms of $\omega$-vectors.   These results extend some of the
 results of \cite[Section 5]{linking}.

  As in Proposition \ref{spike}, let $\lambda_i = a_i(p-1)+b_i$,
   where $a_i \ge 0$, $1 \le b_i \le p-1$.  Following
   \cite{linking}, for $R= ((b_1+1)p^{a_1}-1, \ldots,
   (b_n+1)p^{a_n}-1)$ we call the Milnor basis element $P(R)$
   the {\em Milnor spike} associated to $\lambda$. We note that
   $\omega(R)= \gamma$.  A Milnor spike is an admissible monomial
   \cite{Monks}.  For   example, if $p=3$ and $\lambda = (4,3,1)$
   then the corresponding Milnor spike is $P(8,5,1) =
   P^{32}P^8P^1$,  and for the ${\bf T}$-conjugate partition
   $\gamma = (5,3)$ it is $P(17,5) = P^{32}P^5$. In this example,
   $\lambda_{(1)}' = (3,2)$ and $\lambda_{(2)}' = (2,1)$, so that
   $p(\lambda') = w(3)w(2)\cdot (w(2)w(1))^3 = [x_1, x_2^3,
   x_3^9]\cdot [x_1, x_2^3]^4 \cdot x_1^3$.

\begin{theorem} \label{Milnor}
Let $\lambda$ be ${\bf T}$-regular  with ${\bf T}$-conjugate
$\gamma$.
\begin{itemize}
\item[\rm(i)] $P(R)s(\lambda)= (-1)^{\epsilon(\lambda)}P(R)v(\lambda)
= p(\lambda')$, where $P(R)$ is the Milnor spike associated to
$(\lambda_2, \ldots, \lambda_n)$.
\item[\rm(ii)] If $\lambda= (p-1)\kappa$, where $\kappa$ is column
$2$-regular, $P(S)w(\lambda') = p(\lambda')$, where $P(S)$ is the
Milnor spike associated to $(\gamma_2, \ldots, \gamma_m)$.

\item[\rm(iii)] There are formulae corresponding to {\rm (i)} and
{\rm (ii)} for the  Milnor spikes associated to $\lambda$ and
$\gamma$, with $p(\lambda')$ replaced by  $p(\lambda')^p$.
\end{itemize}
\end{theorem}

\begin{remark} \label{Milnorremark} {\rm
{\rm (iii)} follows immediately from {\rm (i)} and {\rm (ii)} for
degree reasons. The omission of the first terms in $R$ and $S$
corresponds to omitting the highest Steenrod power $P^d$ in the
admissible monomial forms of $P(R)$ and $P(S)$. In fact $d= {\rm
deg}\, p(\lambda')$, so that $P^d p(\lambda')= p(\lambda')^p$. In
the example $p=3$, $\lambda= (4,3,1)$ above, {\rm (i)} states
that $P^8P^1(x_1^8x_2^5x_3) = -P^8P^1( x_1^2\cdot
[x_1,x_2^3]^2\cdot [x_2,x_3^3]) = [x_1,x_2^3,x_3^9]\cdot
[x_1,x_2^3]^4 \cdot x_1^3$.  The case $\lambda= (4,3,1)$ is
excluded from {\rm (ii)}, but in fact $P^5 w(\lambda') =
-p(\lambda')$. We believe that {\rm (ii)} holds, up to sign, for
all ${\bf T}$-regular $\lambda$. }
\end{remark}

 We begin by proving the equivalence of the two statements in {\bf
 (i)}.  For this  we use the following generalization of
 \cite[Theorem 5.9(i)]{linking}.  The proof is based on Lemma
 \ref{Cartan3}, and follows that given in \cite{linking}.
\begin{theorem} \label{omegai}  Let  $R=(r_1, \ldots ,r_t)$ and let
$\omega(R) =\rho$. If the $\omega$-vector $\sigma$ of $x_1^{s_1} \cdots x_n^{s_n}$  does not dominate $\rho$, then
$P(R)(x_1^{s_1} \cdots x_n^{s_n}) = 0$. \qed
\end{theorem}
{\bf Proof of Theorem \ref{Milnor}(i)}\qua By Proposition
\ref{Weylmod1}, if the  monomial $f$  occurs in $v(\lambda)$ and
$f \neq s(\lambda)$, then $\omega(f) \prec \gamma$.   If $R= (r_1,
\ldots, r_n)$ where $r_i= (b_i+1)p^{a_i}-1$, so that $P(R)$  is
the  Milnor spike associated to $\lambda$, then, as noted above,
$\omega(R)= \gamma$. Hence, by Theorem \ref{omegai}, $P(R)$ takes
the same value on $v(\lambda)$ and on its leading term
$(-1)^{\epsilon(\lambda)}s(\lambda)$.

We evaluate  $P(R)s(\lambda)$ by induction on the length $m$ of
$\gamma$.  The base case $m=1$ holds by our previous results, as
follows.  In this case, $\lambda = (p-1, \ldots, p-1, b)$, with $1
\le b \le p-1$,  and has length $n$, while {\rm (i)} states that
$P(R)s(\lambda) = w(\lambda')$, where $R=(p-1, \ldots, p-1, b)$
has length $n-1$.  By Proposition \ref{hatrel}{\rm (ii)}, $P(R)g =
\Hq{(b+1)p_{n-1}-(n-1)}g$ when $\deg g \le (n-1)(p-1)+b$, and we
may choose $g=s(\lambda)$.  Hence the result follows from Theorem
\ref{Minhtrick3}.

For the inductive step, we use  Proposition \ref{spike}{\rm (ii)}
to write  $s(\lambda) = f^p\cdot g$, where $g =v(\lambda_{(1)}) $
and $f = s(\lambda^-)$.  Hence $P(R)s(\lambda)= \sum (P(S)f)^p
\cdot P(T)g$ by Lemma \ref{Cartan3}, where the sum is over
sequences $S=(s_2, \ldots, s_n)$, $T=(t_2, \ldots, t_n)$ such that
$r_i =ps_i +t_i$ for $2 \le i \le n$. Thus $t_n = b_1$, $s_n=0$
and   $t_i \ge p-1$ for $2 \le  i \le n-1$.  If $t_i \neq p-1$ for
some $i<n$, then $P(T)$ has excess $\sum_i t_i > \deg
v(\lambda_{(1)})= \gamma_1$, so that $P(T)(v(\lambda_{(1)}))=0$.
Hence we may assume that $T = (p-1, \ldots, p-1, b_1)$, so that
$s_i = (b_i+1)p^{a_i-1} -1$ for $2 \le  i \le n-1$.  By the
argument for the case $m=1$, $P(T)(v(\lambda_{(1)}))=
w(\lambda'_{(1)})$, and by the induction hypothesis applied to
$\lambda^-$, $P(S)s(\lambda^-) = p(\lambda^-)$.  Since $p(\lambda)
= w(\lambda'_{(1)})\cdot  p(\lambda^-)^p$, the induction is
complete.  \endproof

{\bf Proof of Theorem \ref{Milnor}{\bf (ii)}}\qua  Let $\lambda=
(p-1)\kappa$, where $\kappa$ is   column 2-regular. Then
$\gamma=(p-1)\kappa'$ has length $m=\kappa_1$, and $\lambda_{(i)}
= ((p-1)^{\kappa_i'})$, so that $w(\lambda_{(i)}) =
w(\kappa_i')^{p-1}$.  Also  $S =(p^{\kappa'_2}-1, \ldots,
p^{\kappa'_m}-1 )$, so that $P(S) = P^{t_2}\cdots P^{t_m}$, where
$t_m =  p^{\kappa'_m}-1$ and $t_i = pt_{i+1} + p^{\kappa'_i}-1$
for $1 \le i <m$. We shall argue by induction on $m$, the case
$m=1$, where $P(S)=1$, being trivial.  For $2 \le i \le m$,  let
 $$
  W_i(\lambda') = w(\lambda'_{(1)}) \cdots  w(\lambda'_{(i)})
  \cdot w(\lambda'_{(i+1)})^p \cdots w(\lambda'_{(m)})^{p^{m-i}},
 $$
 so that  $W_1(\lambda') = p(\lambda')$ and $W_m(\lambda') =
w(\lambda')$. We assume as inductive hypothesis on $j$ that
$P^{t_j}W_j(\lambda') =W_{j-1}(\lambda')$ for $j>i$, and prove
this for $j=i$.

It follows from Lemma \ref{Minhlemma} that $P^r(w(n)^{p^i}) = 0$
unless $r= p^i(p_n-p_j)$, where $0\le j \le n$.  The largest of
these values, equal to the degree of $w(n)^{p^i}$,  is $p^i\cdot
p_n$.

Since  $ w(\lambda'_{(i)})$ has degree $p^{\kappa'_i}-1$, it
follows by (downward) induction on $i$ that $t_i$ is the  degree
of $w(\lambda'_{(i)}) \cdot w(\lambda'_{(i+1)})^p \cdots
w(\lambda'_{(m)})^{p^{m-i}}$. We may express $t_i$ explicitly as
the sum
\begin{equation} \label{sol1}
t_i = \sum_{k=i}^m  p^{k-i}(p^{\kappa'_k}-1).
\end{equation}
Hence one term in the expansion of $P^{t_i}(W_i(\lambda'))$ using
the Cartan formula is $W_{i-1}(\lambda')$. We shall complete the
proof by using Lemma \ref{Minhlemma} to show that all  other
terms  in the expansion vanish.  Thus we have to consider the
possible ways to write $t_i$ so that
\begin{equation} \label{sol2}
(p-1)t_i = \sum_{v=1}^{p-1}\left(\sum_{k=1}^{i-1}
(p^{\kappa'_k}-p^{j_{k,v}}) + \sum_{k=i}^m
p^{k-i}(p^{\kappa'_k}-p^{j_{k,v}})\right)
\end{equation}
where $0 \le j_{k,v} \le \kappa'_k$ for $1 \le k \le m$. Equating
(\ref{sol1}) and (\ref{sol2}) and simplifying, we obtain
\begin{equation} \label{sol3}
(p-1)\left(\sum_{k=1}^{i-1} p^{\kappa'_k} + \sum_{k=i}^m
p^{k-i}\right) = \sum_{v=1}^{p-1}\left( \sum_{k=1}^{i-1}
p^{j_{k,v}} +  \sum_{k=i}^m p^{k-i}\cdot p^{j_{k,v}}\right).
\end{equation}
Since $\kappa$ is column 2-regular, $\kappa'$ is strictly
decreasing and so  $\kappa'_{i-1} > \kappa'_i \ge \kappa'_m+m-i
>m-i$.  Hence the $m$ powers of $p$ occurring in the left side of
(\ref{sol3}) are  distinct. By uniqueness of base $p$ expansions,
there are also $m$ distinct powers  on the right  of (\ref{sol3})
and these are a permutation of the powers  on the left.  The
argument is now completed as in the case $p=2$ \cite[Section
5]{linking}. \endproof

We end with evaluations of certain Milnor basis elements on
monomials. While \cite[Lemma 5.6]{linking} generalizes easily to
odd primes, this does not seem to be so useful here as the
following  (weak) generalization of \cite[Proposition
5.8]{linking}.
\begin{proposition} \label{basecase}
Let $R = (r_1, r_2, \ldots)$ where $r_i =p-1$ if $i=b_1, \ldots,
b_m$ and $r_i =0$ otherwise.  Then
 $$P(R)(x_1 \cdots x_n)^{p-1} = 
\begin{cases}
{[x_1^{p^{b_1}}, \ldots, x_n^{p^{b_n}}]^{p-1}} & \text{if $m=n$}, \\
{[x_1, x_2^{p^{b_1}}, \ldots, x_n^{p^{b_{n-1}}}]^{p-1}} & \text{if
 $m=n-1$}.
\end{cases}$$
\end{proposition}

{\bf Proof}\qua This is proved by induction on $|R|$.  The base
of the induction is  Theorem \ref{detp-1}, which is the case
$m=n-1$, $b_i =i$ for $1 \le i \le n-1$. Given a sequence $R =
(r_1, \ldots, r_{j-1}, 0, p-1, p-1, \ldots, p-1)$, let $R'= (r_1,
\ldots, r_{j-1}, p-1,0, p-1, \ldots, p-1)$, so that $|R|-|R'| =
(p-1)(p^{j+1}-1)-(p-1)(p^j-1) = (p-1)^2p^j$. We claim that
$P^{p^j(p-1)} \cdot P(R')$ and $P(R)$ have the same value on any
polynomial of degree $n(p-1)$. To prove this, we use Milnor's
product formula to expand $P^{p^j(p-1)} \cdot P(R')$ in the Milnor
basis.  The Milnor matrix $$
\begin{array}{r|cccccccc}
&r_1 & \ldots & r_{j-1} &0&0&p-1 & \ldots & p-1 \\ \hline 0&0 &
\ldots &0 & p-1 & 0 &0 & \ldots & 0
\end{array}
$$ shows that $P(R)$ occurs with coefficient 1 in the
product. Since $P(R)$ is the unique Milnor basis element of
minimal excess $(n-1)(p-1)$ in degree $|R|$, this proves our claim.

Applying the induction hypothesis to $P(R')$, we have $P(R)(x_1
\ldots x_n)^{p-1} = P^{p^j(p-1)} [x_1, x_2^{p^{b_1}}, \ldots,
x_i^{p^j}, \ldots, x_n^{p^{b_{n-1}}}]^{p-1}$ where $R$ and $R'$
differ in the $i$th term, i.e.\ $b_i = j$ for $R'$ and $b_i=j+1$
for $R$.  By the Cartan formula, this is $[x_1, x_2^{p^{b_1}},
\ldots, x_i^{p^{j+1}}, \ldots, x_n^{p^{b_{n-1}}}]^{p-1}$, and this
completes the induction for the case $m=n-1$.  The case $m=n$ is
proved similarly.  \endproof

 Proposition \ref{basecase} serves as the base of induction for
  the following generalization of \cite[Theorem 5.9(ii)]{linking}
  to odd primes.  The proof, by induction on the length of the
  $\omega$-vector $\sigma$, is essentially the same as in
  \cite{linking}.

\begin{theorem} \label{omegaii} Let $R_0=(r_0, r_1, \ldots ,r_t)$,  $R= (r_1, \ldots, r_t)$ and  $f=x_1^{s_1} \cdots
x_n^{s_n}$, where the base $p$ expansion of each term $r_i$ and
exponent $s_j$   contains only the digits $0$ and $p-1$. Assume
that $f$ and $R_0$ have the same $\omega$-vector $\sigma$.  Then $
P(R) f = \prod_{k=1}^m \Delta_k^{p^{k-1}(p-1)}$, where $m$ is the
length of $\sigma$ and $\Delta_k = [x_{i_1}^{p^{j_1}}, \ldots,
x_{i_\kappa}^{p^{j_\kappa}}]$ is the Vandermonde determinant of
order $\kappa  =\sigma_k/(p-1)$ defined by the subsequences $(s_{i_1},
\ldots, s_{i_\kappa})$ of $(s_1, \ldots, s_n)$ and $(r_{j_1},
\ldots, r_{j_\kappa})$ of $R_0$ consisting of the terms whose
$k$th base $p$ place is $p-1$.
\end{theorem}

\begin{example} \label{omegaiiex} {\rm
Using the  tables
 $$\begin{array}{ccc}
\begin{array}{c|ccc} r_0&p-1&0&p-1\\ \hline r_1&p-1\\ r_2&p-1\\ \hline \sigma &3(p-1)&0&p-1
\end{array}
&\hspace{-5pt}
\begin{array}{c|ccc} x_1&p-1&0&p-1\\  x_2&p-1\\ x_3&p-1\\ \hline \sigma &3(p-1)&0&p-1
\end{array}
&\hspace{-5pt}
\begin{array}{c|ccc} r_0&p-1\\ \hline r_1&p-1&0&p-1\\ r_2&p-1\\ \hline \sigma &3(p-1)&0&p-1
\end{array}  \end{array}
$$
 we obtain $P(p-1,p-1)x_1^{(p^2+1)(p-1)}x_2^{p-1}x_3^{p-1}  =
 [x_1,x_2^p,x_3^{p^2}]^{p-1}\cdot x_1^{p^2(p-1)}$ and
 $P((p^2+1)(p-1),p-1)x_1^{(p^2+1)(p-1)}x_2^{p-1}x_3^{p-1} =
 [x_1,x_2^p,x_3^{p^2}]^{p-1}\cdot (x_1^p)^{p^2(p-1)}$. }
\end{example}

\Addresses\recd

\end{document}

From bob@zeus.math.univ-paris13.fr Fri Jul 19 18:18 BST 2002
Date: Fri, 19 Jul 2002 14:34:47 +0200 (MEST)
From: bob@math.univ-paris13.fr
Subject: Minh-Walker final paper

Dear Colin,

Here's the tex file of the Minh-Walker paper.  I've compiled it and
it looks good.  

One thing puzzled me.  It compiles fine on the PC I'm using, but not on
the machine at P13.  I know that has something to do with an earlier
version of amslatex being used at P13, but I couldn't tell what was 
causing the problem.  (I've learned that the new version requires 
\tag{1} rather than \tag1, and that there are slight differences in
vertical spacing, but those are the only differences I've found before
now.)

Bob

==============